\documentclass[letterpaper]{mc2023}
\usepackage{tabls}
\usepackage{cites}
\usepackage{epsf}
\usepackage{appendix}
\usepackage{ragged2e}
\usepackage[top=1in, bottom=1in, left=1in, right=1in]{geometry}
\usepackage{enumitem}
\setlist[itemize]{leftmargin=*}
\usepackage{caption,color}
\usepackage{amsmath, amssymb, epsfig}
\usepackage{subcaption}
\title{A Polynomial Chaos Approach for Uncertainty Quantification of \\Monte Carlo Transport Codes}
\author{%
  \textbf{G.~Geraci$^1$, K.B.~Clements$^{1,2}$, and A.J.~Olson$^1$}\vspace{3pt} \\
  $^1$Sandia National Laboratories  \\
  $^2$Oregon State University  \\ 
   \url{ggeraci@sandia.gov}, \url{clemekay@oregonstate.edu}, \url{aolson@sandia.gov} 
}
\newcommand{\authorHead}{Gianluca Geraci, Kayla B. Clements and Aaron J. Olson}
\newcommand{\shortTitle}{Efficient PC construction for MC RT}
\newcommand{\Var}[1]{\mathbb{V}ar\left[#1\right]}
\newcommand{\EE}[1]{\mathbb{E}\left[#1\right]}

\newcommand{\Nxi}{N_{\xi}}
\newcommand{\Neta}{N_{\eta}}

\newcommand{\Ctot}{\mathcal{C}_{tot}}
\newcommand{\Cxi}{\mathcal{C}_{\xi}}
\newcommand{\Ceta}{\mathcal{C}_{\eta}}

\newcommand{\Sigsqeta}{\sigma_\eta^2}
\newcommand{\Qpoll}{\tilde{Q}_{N_\eta}}

\newcommand{\EExi}[1]{\mathbb{E}_\xi\left[#1\right]}
\newcommand{\EEeta}[1]{\mathbb{E}_\eta\left[#1\right]}
\newcommand{\Vxi}[1]{\mathbb{V}ar_\xi\left[#1\right]}

\newcommand{\Stm}{\Sigma_{t,m}^{0}}
\newcommand{\Std}{\Sigma_{t,m}^{\Delta}}

\begin{document}
\maketitle
\justify 
\parskip 6pt plus 1 pt minus 1 pt

\begin{abstract}
 In this contribution, we discuss the construction of Polynomial Chaos (PC) surrogates for Monte Carlo (MC)
 radiation transport (RT) applications via non-intrusive spectral projection. This contribution focuses on improvements with respect to the approach that we previously introduced in~\cite{GeraciMC2021}. We focus on understanding the impact of re-sampling cost on the algorithm performance and provide algorithm refinements, which allow to obtain unbiased estimators for the variance, estimate the PC variability due to limited samples, and adapt the expansion. An attenuation-only test case is provided to illustrate and discuss the results. 
\end{abstract}
\vspace{6pt}

\section{INTRODUCTION} 
\label{SEC:intro}
Uncertainty Quantification (UQ) deals with the characterization and quantification of the sources of uncertainty that affect the output of a numerical simulation. To compute statistics for the output of interest of a system, a computer code must be evaluated multiple times. Many applications, including RT applications, are challenging for UQ analyses due to the high cost of each computer simulation, which limits the attainable accuracy of the sought-after statistics. In these cases, the construction of a surrogate has the promise of alleviating the overall cost of UQ by orders of magnitude. A surrogate is a mathematical model that, once trained, can be effectively used to replace the original and expensive computational model to obtain a much larger set of evaluations at virtually no additional cost. However, the cost of surrogate construction is often dominated by the construction of its training dataset and dealing efficiently with this is the key to retain computational efficiency. While several surrogates exist, the most well-known in the context of UQ is the so-called Polynomial Chaos (PC)~\cite{LeMaitre2010}, which presents features that are particularly well-suited for UQ (see next section). In this work, we are interested in dealing with the PC construction from the perspective of MC RT solvers. These solvers, like many other non-deterministic solvers encountered in applications like plasma physics, computer networks/cybersecurity, and turbulent flows, have received less attention than their deterministic counterparts in the context of UQ and surrogate construction. MC RT solvers pose a major challenge to standard existing approaches: the solvers' variability needs to be controlled in order to avoid corrupting the UQ statistics. In this contribution, we build on our previous work~\cite{GeraciMC2021} in which we discussed how PC can be efficiently built in the presence of limited number of particle histories for each MC RT simulation. In~\cite{GeraciMC2021}, we demonstrated that the use of a single particle history, for each MC RT simulation, provably leads to a smaller variance than over-resolving each single simulation, for the same number of total particle histories. Here, we extend~\cite{GeraciMC2021} by several contributions. First, we generalize the above reported result by accounting for a re-sampling cost corresponding to the generation of a new system's configuration, which depends on the uncertain parameters. Second, we derive an approach to obtain the surrogate variability as function of the coefficients' variability due to the sampling process. Notably, this expression is derived by introducing a correction for the covariance term among coefficients that extends our variance deconvolution approach; see~\cite{ClementsCSRI2021,ClementsANS2022}. Third, we propose two algorithmic refinements to the previous approach that allow to obtain unbiased statistics, e.g., for the variance, and automatically select the polynomial expansion to minimize the excitation of spectral terms introduced by, e.g., the under-resolved MC RT simulations used to train the surrogate. Several numerical results are reported to corroborate our findings.

\section{BACKGROUND ON POLYNOMIAL CHAOS} 
\label{SEC:theory}
The PC surrogate is the most used surrogate in UQ due to some of its properties that simplify the UQ analyses. Given a generic Quantity of Interest (QoI) $Q \in \mathbb{R}$ of a system and a vector of uncertain parameters $\xi \in \Xi \subset \mathbb{R}^d$ with joint probability density function $p(\xi)$, PC approximates the QoI via a truncated polynomial series 
\begin{equation}
 Q(\xi) = \sum_{k=0}^\infty \beta_k \Psi_k(\xi) \approx \sum_{k=0}^P \beta_k \Psi_k(\xi), 
\end{equation}
where $\Psi_k$ represents the $k^{\text{th}}$ term of a polynomial basis orthogonal to the measure $p(\xi)$. The number of $P+1$ terms for the series can be determined in different ways, one of which is to prescribe the total polynomial degree $n_0$ of the expansion; e.g., for a full tensor expansion PC has $P+1 = (n_0+d)! / (n_0!\,d!)$ terms. In our previous work~\cite{GeraciMC2021}, we focused on the so-called non-intrusive spectral projection (NISP) approach~\cite{LeMaitre2010} to evaluate these coefficients 
\begin{equation}
  \label{eq:NISP}
  \begin{split}
  \beta_k = \frac{ \EExi{ Q(\xi) \Psi_k(\xi) } }{ b_k }, \quad \mathrm{where} \quad \EExi{\left( \cdot \right)} = \int_{\Xi} \left( \cdot \right) p(\xi) \mathrm{d}\xi,
  \end{split}
\end{equation}
and the normalization term $b_k= \EExi{\Psi_k^2(\xi)}$ is a constant obtained in closed form, once the basis is prescribed. Hereinafter, we use the subscript to define the integration variable, whereas the absence of subscript indicates the integration with respect to all variables of the integrand. PC is widely adopted in UQ because, once the coefficients are obtained, moments (like expected value, variance, \emph{etc.}) and other quantities (like Sobol' indices~\cite{Crestaux2009}) can be directly obtained from them. For instance, the mean and variance of the QoI $Q$ can be estimated as
\begin{equation}
\label{eq:PCEstats}
\begin{split}
 \EE{Q} \approx \beta_0 \quad \mathrm{and} \quad \Var{Q} \approx \sum_{k=1}^P \beta_k^2 \,b_k.
\end{split}
\end{equation}
In the MC RT context, if we consider an elementary event $f$, i.e., an event corresponding to a single MC particle history, a transport solution depends not only on the random vector $\xi$, but also on the vector of variables representing the solver's randomness that we indicate with $\eta$. It follows that $Q$, by definition, corresponds to a statistic of $f$ for each realization of $\xi$, i.e., $Q(\xi) = \EEeta{f}$, such that Eq.~\eqref{eq:NISP} can be manipulated as
\begin{equation}\label{eq:NISP-with-eta}
 \beta_k = \frac{\EExi{ Q(\xi) \Psi_k(\xi) }}{b_k}  = \frac{\EExi{ \EEeta{f(\xi,\eta)}\Psi_k(\xi) }}{b_k}   = \frac{\EE{f(\xi,\eta) \Psi_k(\xi)}}{b_k},
\end{equation} 
which illustrates the need for obtaining events $f$ in the combined space $\xi-\eta$. We note here that the use of the variable $\eta$ is notional; this variable is indeed unknown and only represents all the solver's randomness in the MC RT simulation. Moreover, its dimensionality can be arbitrarily large without affecting the following results because it would never need to be sampled. In~\cite{GeraciMC2021}, we discussed how the previous expression can be sampled to obtain a computable expression for the PC coefficients. Two sampling strategies, either \emph{nested} or \emph{direct}, can be obtained from Eq.~\eqref{eq:NISP-with-eta} as
%
\begin{equation}
\label{eq:NISP_est}
\beta_k 
    \approx \frac{1}{b_k} \frac{1}{N_\xi} \sum_{i=1}^{N_\xi} \left( \Psi_k(\xi^{(i)}) Q^{(i)} \right)
    \triangleq \hat{\beta}_k
     \approx \frac{1}{b_k} \frac{1}{N_\xi} \sum_{i=1}^{N_\xi} \left( \frac{\Psi_k(\xi^{(i)})}{N_\eta} \sum_{j=1}^{N_\eta} f(\xi^{(i)},\eta^{(j)}) \right) 
     \triangleq \hat{\beta}_{k,N_\eta}, 
\end{equation}
with $\hat{\beta}_k$ being the NISP on the unpolluted QoI and $\hat{\beta}_{k,N_\eta}$ being its approximation using the MC RT quantities, which, in turn, depend on the choice of $\Neta$; the direct approach is obtained for $\Neta=1$. As a shorthand, we also define $\Qpoll = (\sum_{j=1}^{N_\eta} f(\xi^{(i)},\eta^{(j)})) / \Neta$ as the polluted approximation of the QoI with $\Neta$ particles.

\subsection{Global Sensitivity Analysis}
PC is a popular surrogate because it allows for the effortless evaluation of several quantities directly from its coefficients. A prominent example is represented by the evaluation of Sobol' indices~\cite{Sobol,Saltellibook}, which are used to measure the sensitivity of a model's output with respect to its input variables. This task is referred to as global sensitivity analysis (GSA). Given a group of variables $u$, sensitivity indices can be defined~\cite{Crestaux2009} as 
\begin{equation}\label{eq:pce-si}
 S_u = \frac{\sum_{k \in K_u} \beta_{k}^2 \, b_k}{ \sum_{k=1}^P \beta_{k}^2 \, b_k } \quad \mathrm{and} \quad S_{T_i} = \sum_{ u \ni i  } S_u,
\end{equation}
where $K_u$ represents all the PC expansion terms in which all the variable(s) in $u$ are present, while the total sensitivity index $S_{T_i}$ includes all first order indices that contain the $i$th variable. High-order sensitivity indices can also be obtained from the PC coefficients, e.g., to quantify the importance of variables with respect to the skewness and kurtosis; see~\cite{Geraci2016}. 

\section{NOVEL CONTRIBUTIONS}
\label{SEC:contributions}
\subsection{Estimator Performance with Re-sampling Cost }
In~\cite{GeraciMC2021}, we demonstrated that for a cost model equal to $\Ctot=\Nxi \times \Neta$, the variance of a coefficient is minimized when $\Neta = 1$ and, as a consequence, $\Nxi = \Ctot$. This result might be difficult to adapt to practical applications in which at each new sample in the parameters' space it corresponds a configuration that needs to be initialized, e.g., by meshing a geometry. In these cases, a \emph{re-sampling cost} would need to be included. We generalize here the cost model as $\Ctot = \Nxi \left( \Cxi + \Ceta \Neta \right)$, where $\Cxi$ measures the re-sampling cost and $\Ceta$ expresses the unitary cost for the evaluation of a single particle in the MC RT code\footnote{The previous model can be recovered for $\Ceta=0$ (no re-sampling cost) and unitary particle's cost, i.e., $\Ceta=1$.}. If both numbers of UQ samples and particles are treated as real numbers, for a fixed cost $\Ctot$, the number of UQ samples is a function of the particles $\Neta$, i.e., $\Nxi = \Ctot / \left( \Cxi + \Ceta \Neta \right)$. From this we can express the variance of the estimator for the $k$th PC coefficient as (from Eq.~\eqref{eq:NISP_est})
\begin{equation}
\label{eq:var_beta_cost}
 \Var{\hat{\beta}_{k,N_\eta}} = \frac{1}{b^2_k} \frac{\Cxi+\Ceta \Neta}{\Ctot} 
                                \left( \Vxi{ Q \Psi_k } + \frac{ \EExi{ \Psi_k^2 \Sigsqeta } }{\Neta} \right).
\end{equation}
We are interested in determining until which conditions the direct approach, corresponding to $\Neta=1$ produces a lower variance than the nested approach, which uses $\Neta=\Neta^{nested}>1$, for a prescribed total cost. After several manipulations, it is possible to show that if the following condition is satisfied
\begin{equation}
 \Neta^{nested} > \frac{\EExi{ \Psi_k^2 \Sigsqeta } }{ \Vxi{ Q \Psi_k } } \frac{\Cxi}{\Ceta},
\end{equation}
then the direct approach is more efficient than the nested one. For instance, if the re-sampling cost is negligible, i.e., $\Cxi=0$, the direct approach is more efficient than the nested one if $\Neta^{nested}>0$, which is always satisfied as we previously demonstrated; see~\cite{GeraciMC2021}.

\subsection{Unbiased Formulation}
\label{SEC:bias_corr}
The first algorithmic refinement that we present consists of an unbiased estimation for the evaluation of statistics like the variance and, as a consequence, all quantities derived from the variance, such as Sobol' indices. Following Eq.~\eqref{eq:PCEstats}, the variance requires the squared coefficients. Since we are employing a sampling approach to obtain the coefficients, simply squaring the estimator in Eq.~\eqref{eq:NISP_est} produces a biased estimator. It can be shown that the bias can be eliminated by subtracting the variance of the estimator, i.e, 
\begin{equation}
 \Var{Q} \approx \sum_{k=1}^P \left( (\hat{\beta}_{k,N_\eta})^2-\Var{\hat{\beta}_{k,N_\eta}} \right) \,b_k,
\end{equation}
where we use $\left( \hat{\beta}_{k,N_\eta} \right)^2-\Var{\hat{\beta}_{k,N_\eta}}$ as an unbiased estimator of $\beta_k^2$, where $\Var{\hat{\beta}_{k,N_\eta}}$ is obtained as in Eq.~\eqref{eq:var_beta_cost}, and where the ratio $(\Cxi+\Ceta \Neta)/(\Ctot)$ corresponds to $1/\Nxi$.

\subsection{PC Expansion Trim}
We discuss here another algorithmic refinement that is useful to ameliorate the numerical performance of the PC approach illustrated above. We have discussed that, given a maximum polynomial order for the expansion, $P+1$ coefficients are given in the expansion. Nevertheless, this expansion could include terms that should not be included in the expansion; for instance, whenever a larger degree for the expansion is selected with respect to the true underlying QoI, high-order frequencies are spuriously introduced unless the method is able to identify and eliminate these terms. Similarly, in the case of a compressible function in which only a small set of non-zero coefficients is present, the method presented above would not be able to eliminate the additional expansion terms. As illustrated above, the variance of the QoI is estimated by summing the contribution associated to each coefficient, therefore, if the expansion is larger than necessary the estimated variance would be larger than the true variance. Of course, this challenge is particularly important in the presence of under-resolved (in term of UQ samples) applications, since for larger numbers of $\Nxi$ these coefficients would converge to zero. In order to mitigate this challenge, we propose to adopt an algorithm that trims the PC basis \textit{a posteriori}. The proposed algorithm works by re-ordering the variance contributions in decreasing order and truncating the series after the estimated variance reaches the value of the sampling variance obtained via the variance deconvolution estimator; see~\cite{ClementsCSRI2021,ClementsANS2022}.    

\subsection{PC Variability}
In this section, we demonstrate how the PC variability can be determined as a function of the variability of the coefficients' estimators. If the QoI $Q$ is expressed in vector form for a fixed realization of $\xi$, namely $\overline{\xi}$, we can write
\begin{equation}
\label{eq:PCE_variance}
 Q(\xi=\overline{\xi}) = \hat{\beta}_0 + \underline{\Psi}^\mathrm{T}(\overline{\xi}) \underline{\hat{\beta}}_{N_\eta}
 \quad \rightarrow \quad
 \Var{ \underline{\Psi}^\mathrm{T}(\overline{\xi}) \underline{\hat{\beta}}_{N_\eta} } =  \underline{\Psi}^\mathrm{T}(\overline{\xi}) \Var{ \underline{\hat{\beta}}_{N_\eta}} \underline{\Psi}(\overline{\xi}), 
\end{equation}
where $\underline{\Psi}(\overline{\xi}) = \left[ \Psi_1(\overline{\xi}), \dots, \Psi_P(\overline{\xi}) \right]^{\mathrm{T}} \in \mathbb{R}^P$ is the vector of basis evaluated at $\overline{\xi}$, $\underline{\hat{\beta}}_{N_\eta} \in \mathbb{R}^{P}$ is the vector of coefficients' estimators (not including the mean term $k=0$), and $\Var{ \underline{\hat{\beta}}_{N_\eta}} \in \mathbb{R}^{P \times P}$ is their variance-covariance matrix. This latter term is due to the fact that each coefficient is evaluated by resorting to the same set of MC RT runs, therefore, in general, their covariance is different from zero. The diagonal of the variance-covariance matrix corresponds to the variance in Eq.~\eqref{eq:var_beta_cost}, whereas the off-diagonal $k,r$ term can be written as
\begin{equation}
\begin{split}
 \mathrm{C}ov\left[ \hat{\beta}_{k,N_\eta}, \hat{\beta}_{r,N_\eta} \right] &= \EE{ \hat{\beta}_{k,N_\eta} \hat{\beta}_{r,N_\eta}} - \EE{ \hat{\beta}_{k,N_\eta} } \EE{ \hat{\beta}_{r,N_\eta} } \\
 &= \frac{1}{\Nxi} \left( \frac{\EE{\Psi_k \Psi_r \Qpoll^2}}{b_k b_r} 
 - \frac{ \EE{ \Psi_k \Qpoll } \EE{ \Psi_r \Qpoll } }{b_k b_r} \right).
\end{split} 
\end{equation}
We now note that this covariance includes the effect of the variability induced by both $\xi$ and $\eta$ and, therefore, the variability induced by $\eta$ is only due to the use of the MC RT code. By using properties for both the mean and variance, we can separate the effect of $\eta$ and write
\begin{equation}
 \mathrm{C}ov\left[ \hat{\beta}_{k,N_\eta}, \hat{\beta}_{r,N_\eta} \right] = 
 \mathrm{C}ov_{\xi}\left[ \hat{\beta}_{k}, \hat{\beta}_{r} \right] + 
 \frac{1}{N_\xi} \frac{1}{b_k b_r} \EExi{ \Psi_k \Psi_r \frac{\sigma_\eta^2}{N_\eta} },
%
\end{equation}
which expresses how the covariance term due to the uncertain parameters is affected by the presence of the MC RT noise. It follows that, at each location $\overline{\xi}$, we can use Eq.~\eqref{eq:PCE_variance} to obtain the variance and standard deviation of the PC response due to the use of random samples in the UQ parameters without the inclusion of the additional noise due to the MC RT code.

\section{NUMERICAL RESULTS} 
\label{SEC:num_res}
In this section, we present numerical results that illustrate some of the novel algorithmic capabilities introduced in this work to extend the approach presented in~\cite{GeraciMC2021}. We consider the stochastic, one-dimensional, neutral-particle, absorption-only, mono-energetic, and steady-state radiation transport equation with a normally incident beam source of magnitude one in a 1D space:
\begin{equation}
\begin{split}
 \mu \frac{\partial \psi(x,\mu,\xi)}{\partial x} + \Sigma_t(x,\xi) \psi(x,\mu,\xi) = 0, \quad \mathrm{where} \quad 0 \leq x \leq L,
\end{split}
\end{equation}
where $\psi$ is the angular flux; $\Sigma_t$ is the total cross section; $x$ and $\mu$ denote spatial and angular dependence; and $\xi$ represents the vector of random parameters. As done in~\cite{GeraciMC2021}, we consider a slab constituted by $d$ material sections with assigned boundaries and an uncertain total cross section for each material described by $\Sigma_{t,m}(\xi_m) = \Stm  + \Std \xi_m$, where $\Stm$ and $\Std$ denote the average value and length of the interval of cross section values for the $m^\text{th}$ material and $\xi_m \sim \mathcal{U}(-1,1)$. We note that our approach can be used as long as the distribution of the parameters $\xi$ is known and a MC RT code can be used to obtain $\Neta$ particle histories for it. As such, we could also consider uncertainties in angular source, magnitude, boundaries, \emph{etc}.; however, this test case is used here due to the presence of a closed-form solution, which can be used for verification. We demonstrate our results for the cases used in~\cite{GeraciMC2021}: for $d=1$ we use $\Delta x_1 = 1$ cm and $\Sigma_{t,1}(\xi_1)\sim\mathcal{U}(0.05, 1.95)$ cm$^{-1}$, whereas for $d=3$ we use $\Sigma_{t,i}(\xi_i)\sim\mathcal{U}(0.01, 0.59)$ cm$^{-1}$ for all material sections. 

To evaluate the statistics, e.g., probability density functions of the estimators, we repeat the numerical experiments for $1,500$ independent repetitions. Moreover, given the availability of the exact solution, we report the mean-squared-error (MSE) defined as $MSE[\hat{X}] = \mathbb{E}[ (\hat{X}-X_{exact})^2 ]$, where $X_{exact}$ is the exact statistics, e.g., the variance, and $\hat{X}$ its approximation obtained either via variance deconvolution or PC. 


\subsection{PC response's variability}
\begin{figure}[hbt!]
     \centering
     \begin{subfigure}[]{0.325\textwidth}
         \centering
         \includegraphics[width=0.99\textwidth]{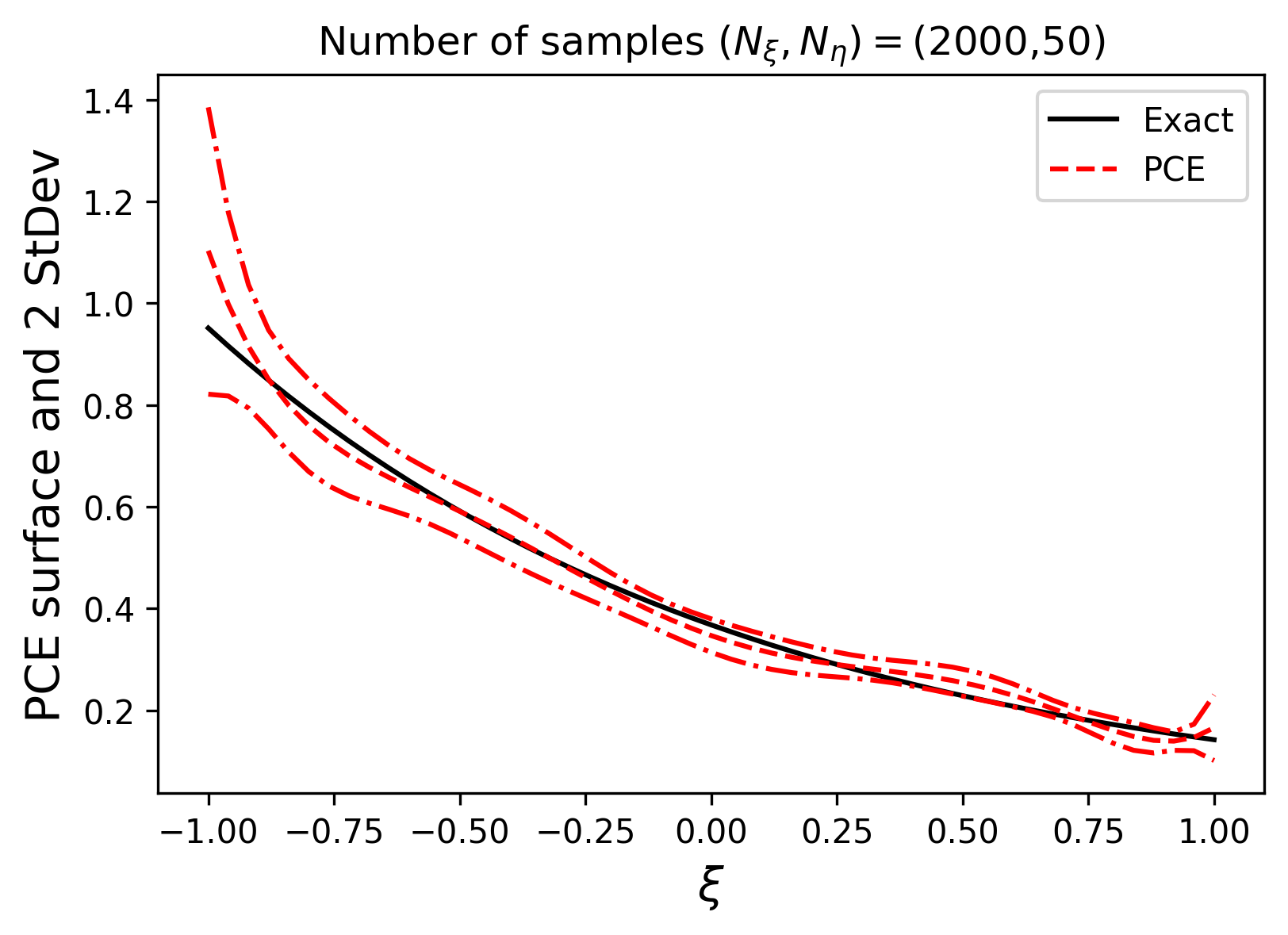}
         \caption{PCE W/O Trim (Sample 1)}
         \label{fig:variance_all}
     \end{subfigure}
     \hfill
     \begin{subfigure}[]{0.325\textwidth}
         \centering
         \includegraphics[width=0.99\textwidth]{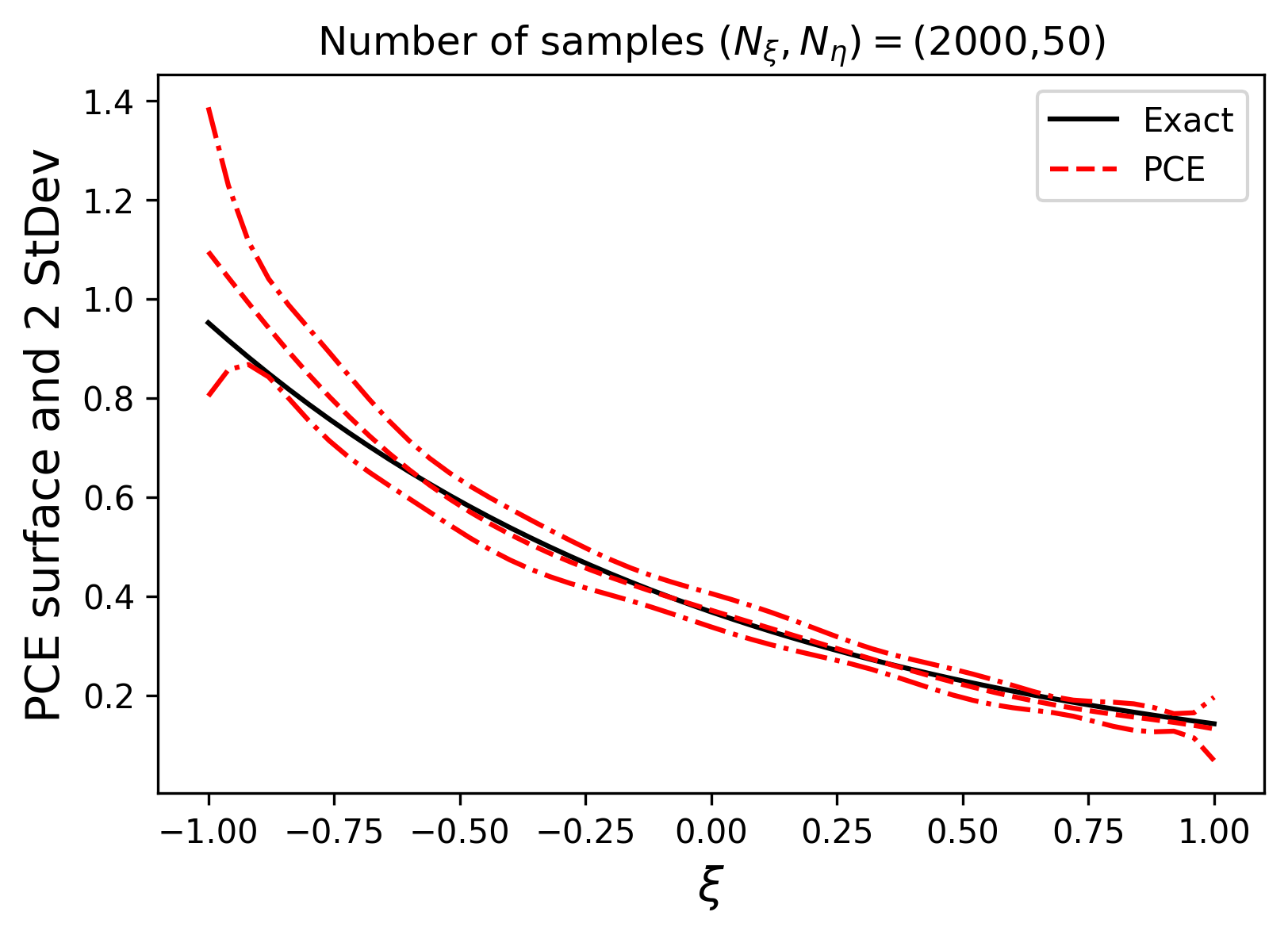}
         \caption{PCE W/O Trim (Sample 2)}
         \label{fig:variance_all}
     \end{subfigure}
     \hfill
     \begin{subfigure}[]{0.325\textwidth}
         \centering
         \includegraphics[width=0.99\textwidth]{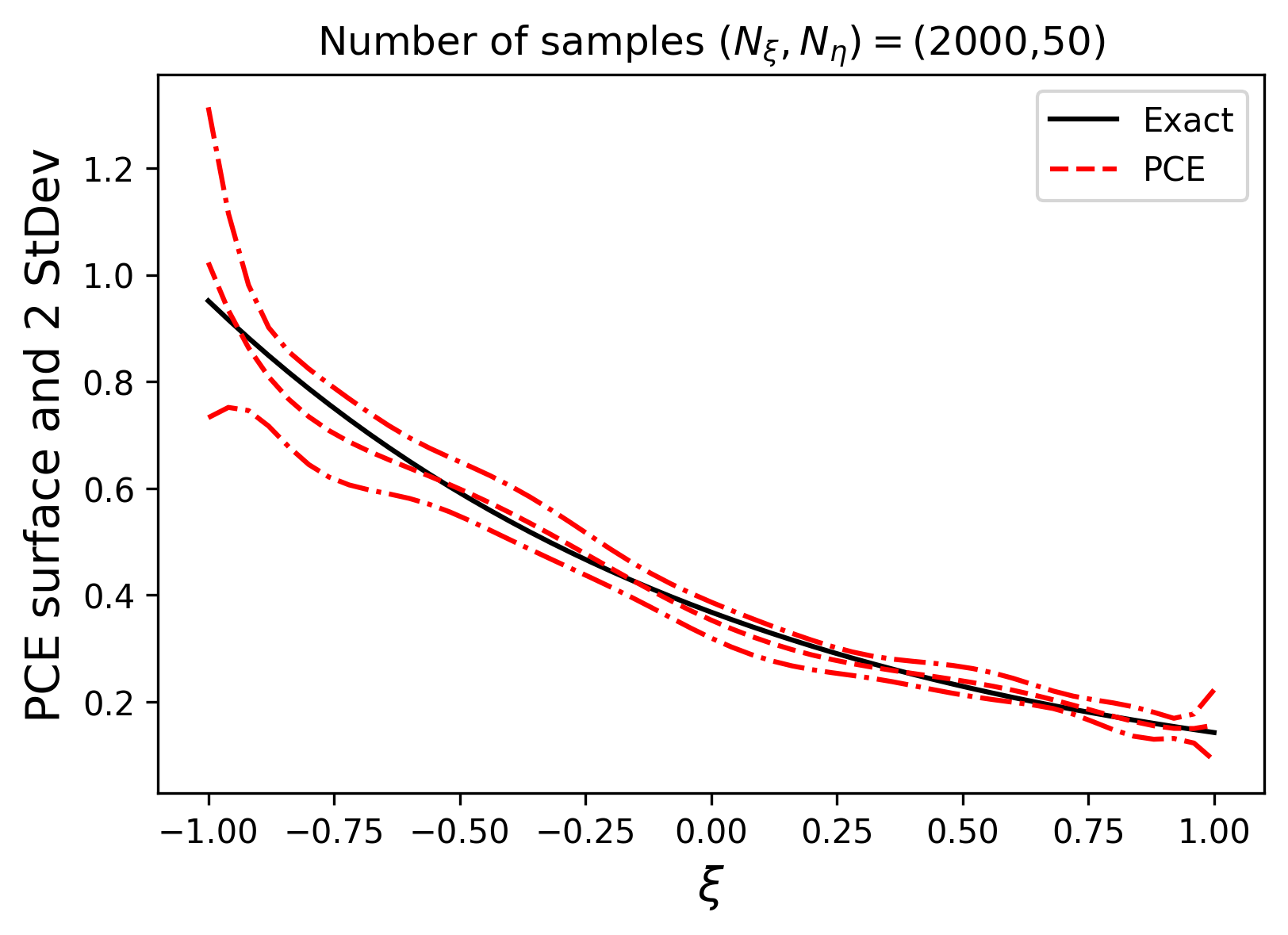}
         \caption{PCE W/O Trim (Sample 3)}
         \label{fig:variance_all}
     \end{subfigure} 
     \begin{subfigure}[]{0.32\textwidth}
     \vspace{1mm}
         \centering
         \includegraphics[width=0.99\textwidth]{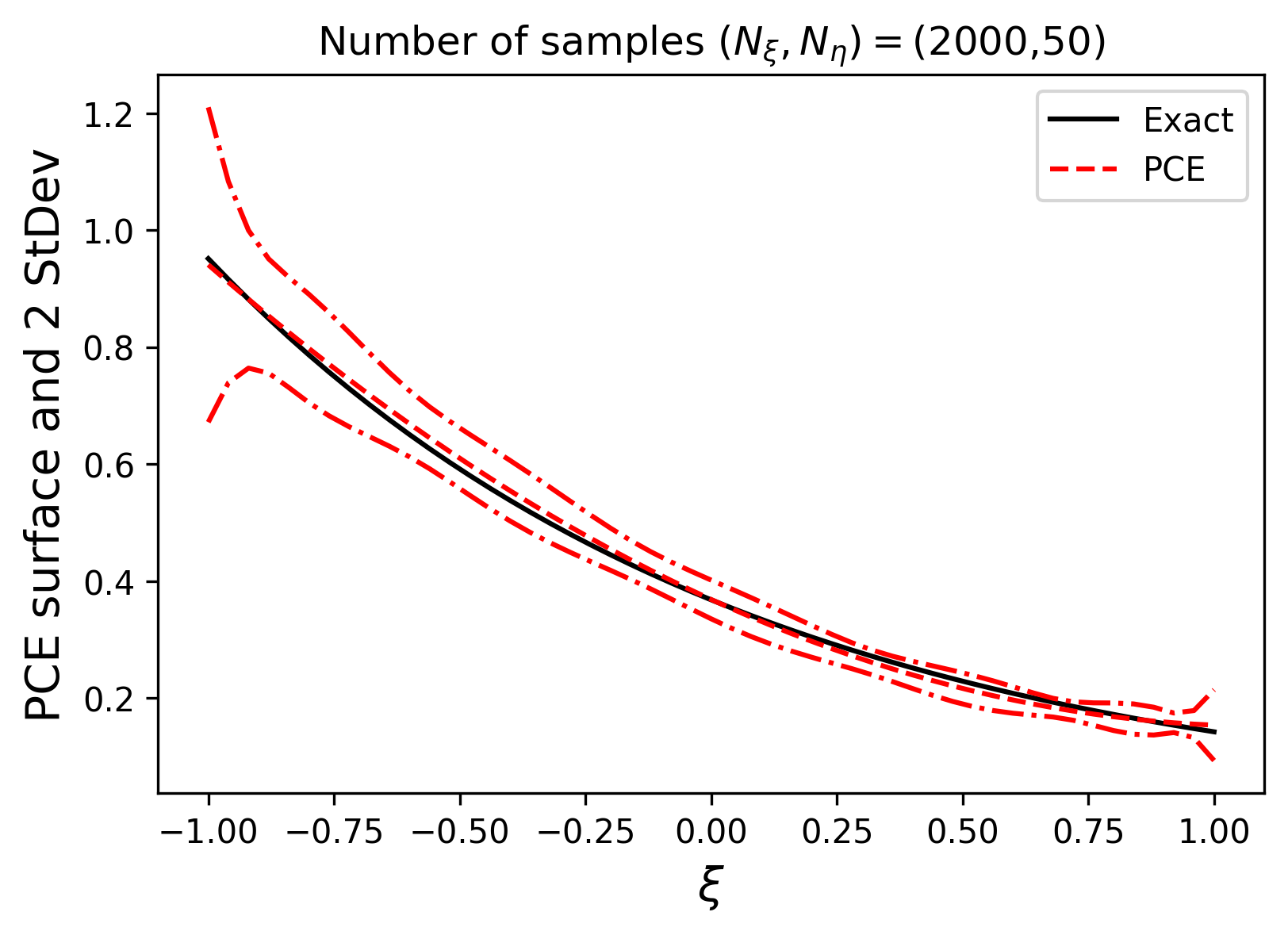}
         \caption{PCE W/ Trim (Sample 1)}
         \label{fig:variance_all}
     \end{subfigure}
     \hfill
     \begin{subfigure}[]{0.32\textwidth}
     \vspace{1mm}
         \centering
         \includegraphics[width=0.99\textwidth]{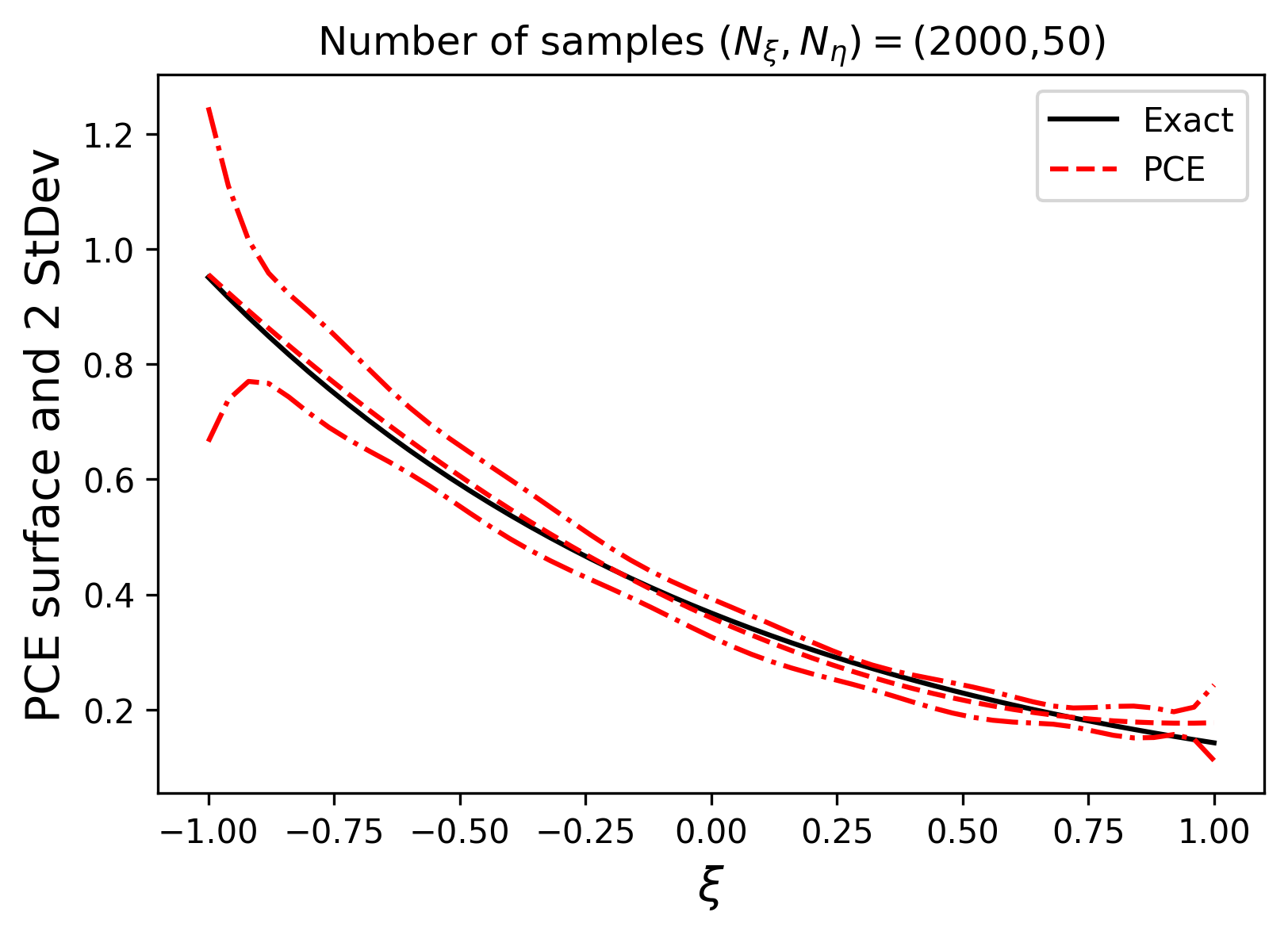}
         \caption{PCE W/ Trim (Sample 2)}
         \label{fig:variance_all}
     \end{subfigure}
     \hfill
     \begin{subfigure}[]{0.32\textwidth}
     \vspace{1mm}
         \centering
         \includegraphics[width=0.99\textwidth]{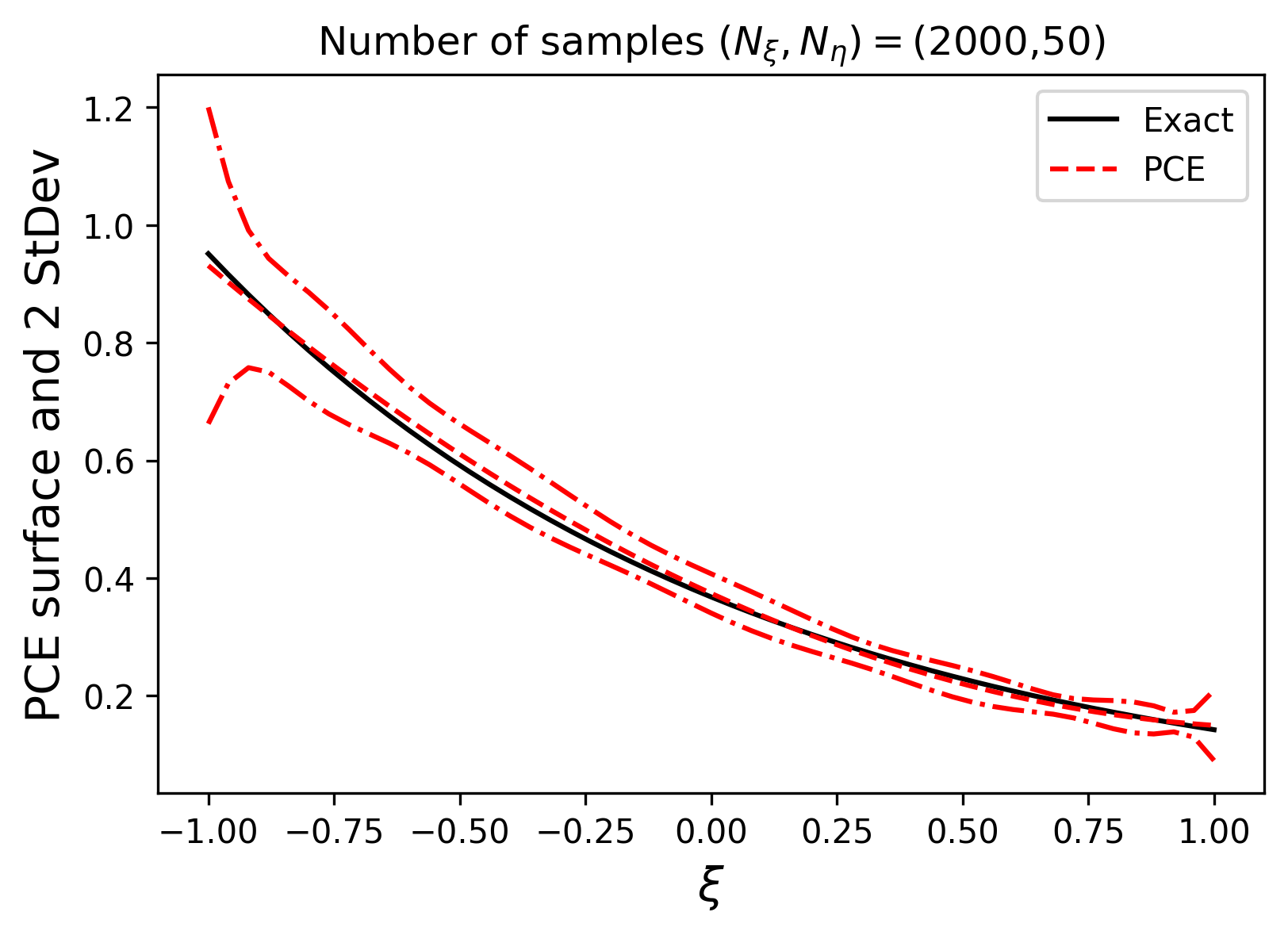}
         \caption{PCE W/ Trim (Sample 3)}
         \label{fig:variance_all}
     \end{subfigure}
     \caption{Three PC repetitions for the 1D attenuation problem (dashed red) with (bottom) and without (top) the expansion trim. All the results are obtained with $\Nxi=2000$ and $\Nxi=50$. The exact attenuation profile is reported in black.}
     \label{fig:1D_responses}
\end{figure}

We start by demonstrating the use of the PC coefficients' covariance to obtain the uncertainty in the PC response. To simplify the presentation, we use the case with a single uncertainty and we report the estimated response of PC for the transmittance along with 2 standard deviations. Moreover, to illustrate the effect of the expansion trim, we compare three repetitions of the algorithm with and without trim with a total polynomial order of 6 for the expansion; see Fig.~\ref{fig:1D_responses}. We note that the expansion trim is able to suppress the spurious oscillations introduced by the high-order expansion terms. 
\subsection{Variance}
\label{SSEC:num_res_variance}
\begin{figure}[hbt!]
     \centering
     \begin{subfigure}[b]{0.32\textwidth}
         \centering
         \includegraphics[width=0.99\textwidth]{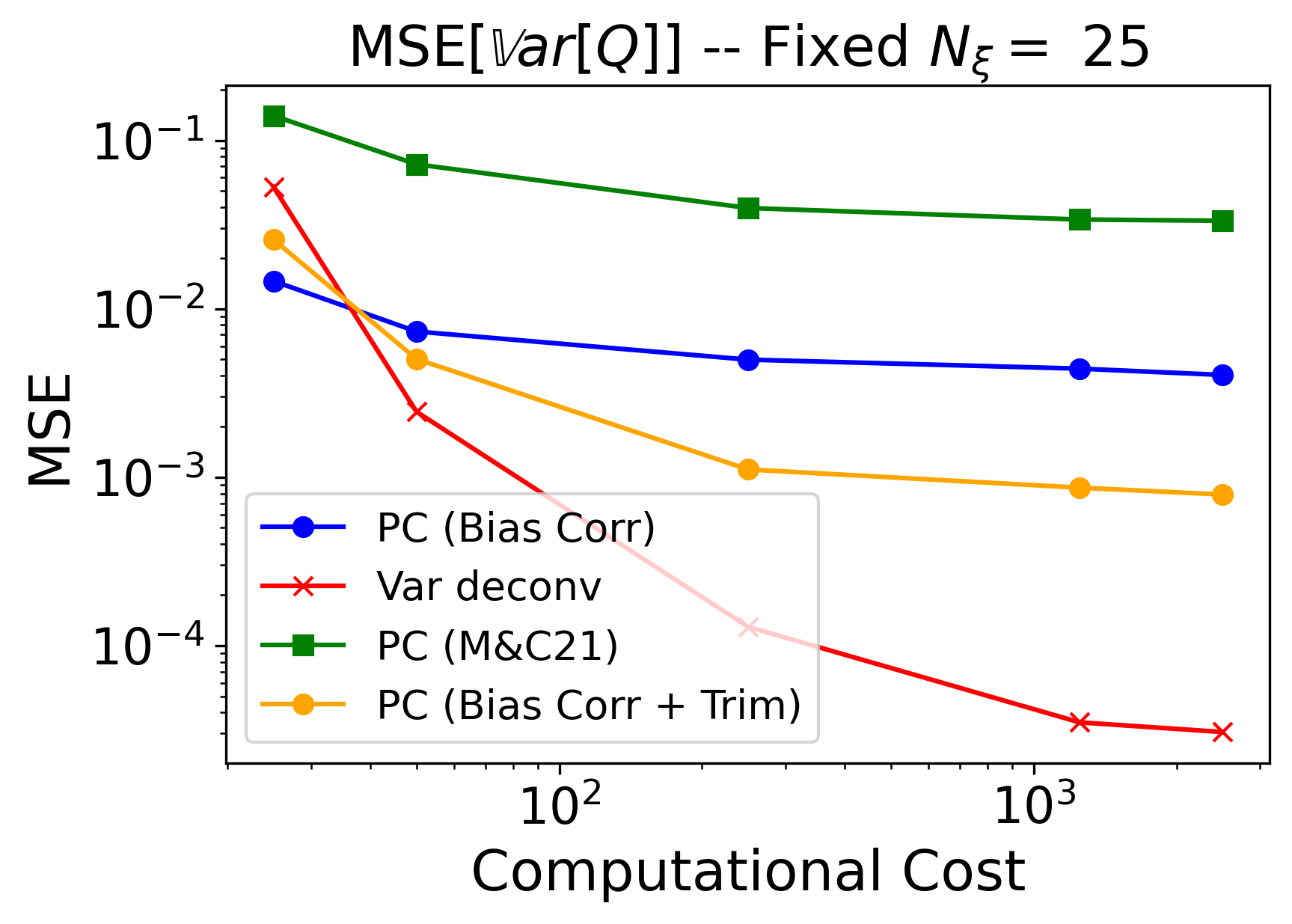}
         \caption{$\Nxi=25$}
         \label{fig:MSE_var_Nxi25}
     \end{subfigure}
     \hfill
     \begin{subfigure}[b]{0.32\textwidth}
         \centering
         \includegraphics[width=0.99\textwidth]{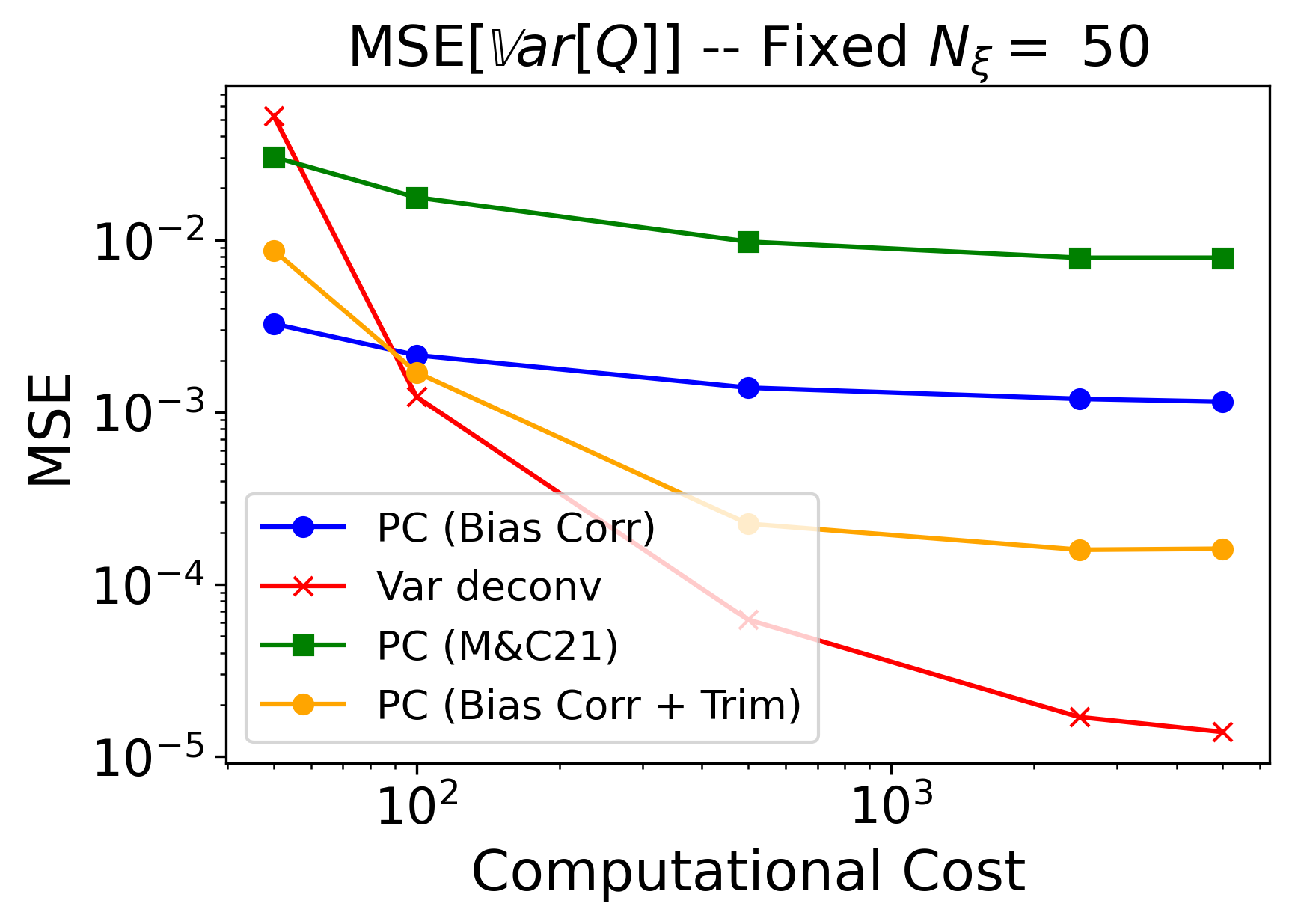}
         \caption{$\Nxi=50$}
         \label{fig:MSE_var_Nxi50}
     \end{subfigure}
     \hfill
          \begin{subfigure}[b]{0.32\textwidth}
         \centering
         \includegraphics[width=0.99\textwidth]{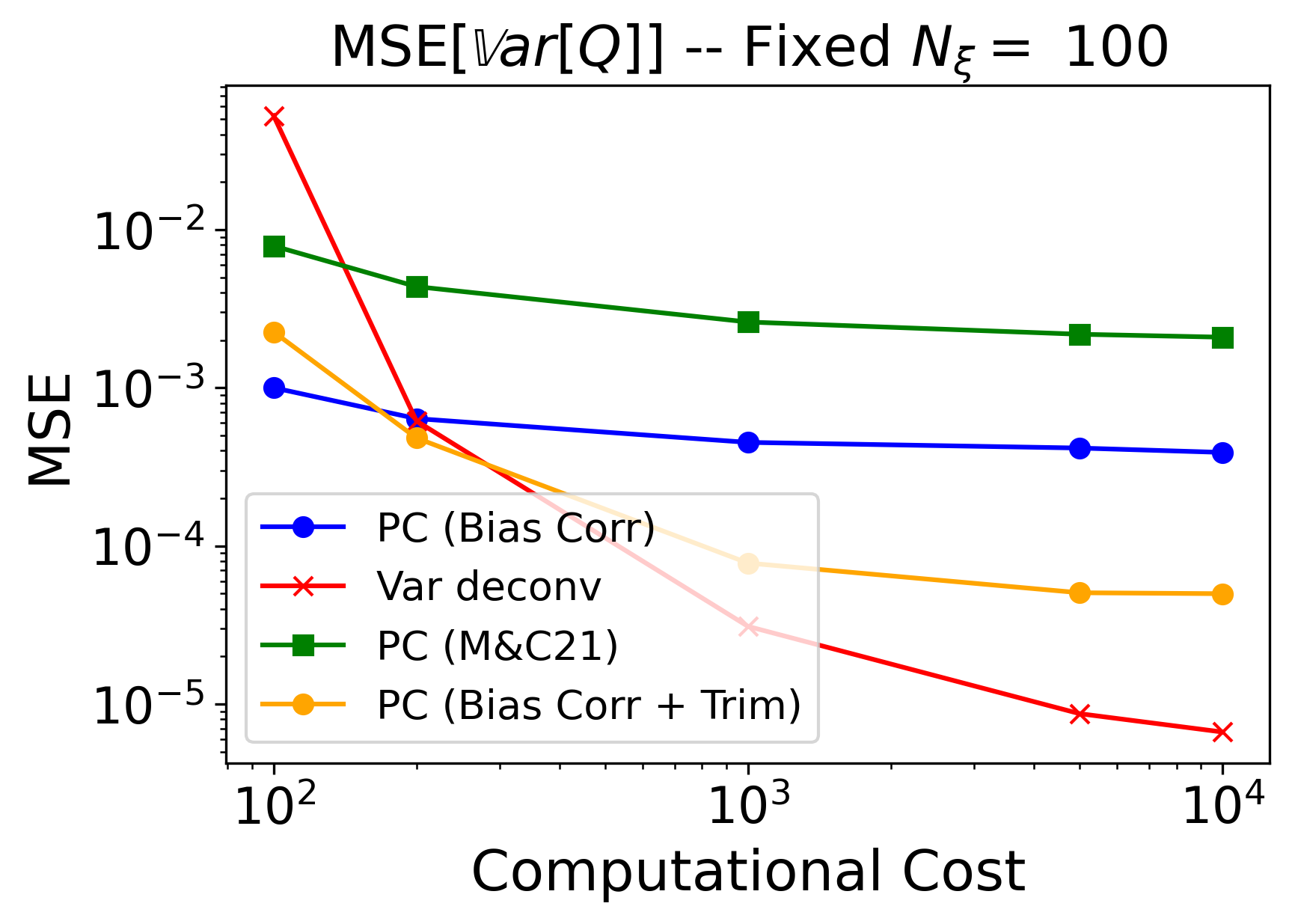}
         \caption{$\Nxi=100$}
         \label{fig:MSE_var_Nxi100}
     \end{subfigure}
     \begin{subfigure}[b]{0.32\textwidth}
         \centering
         \includegraphics[width=0.99\textwidth]{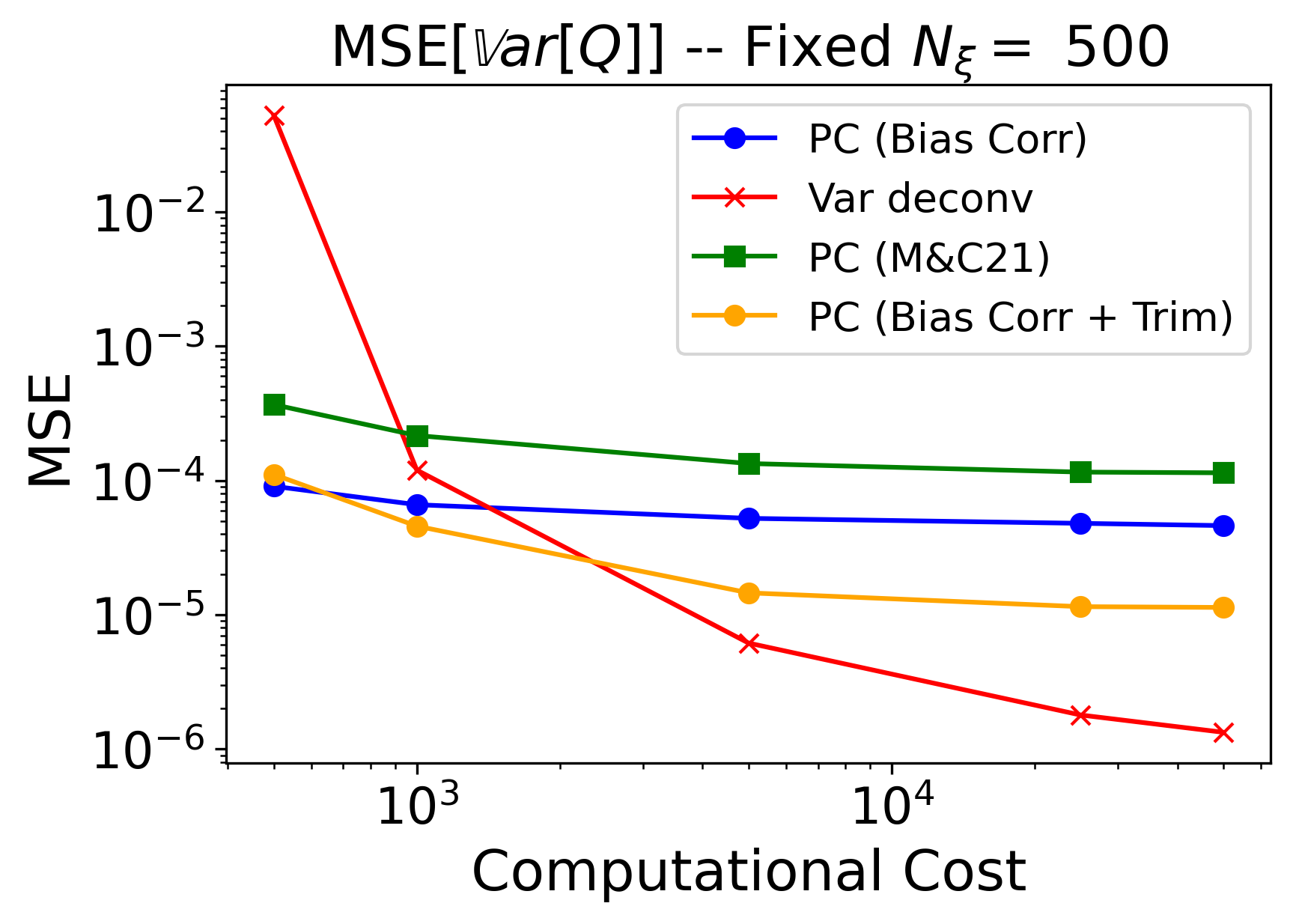}
         \caption{$\Nxi=500$}
         \label{fig:MSE_var_Nxi500}
     \end{subfigure}
     \hfill
          \begin{subfigure}[b]{0.32\textwidth}
         \centering
         \includegraphics[width=0.99\textwidth]{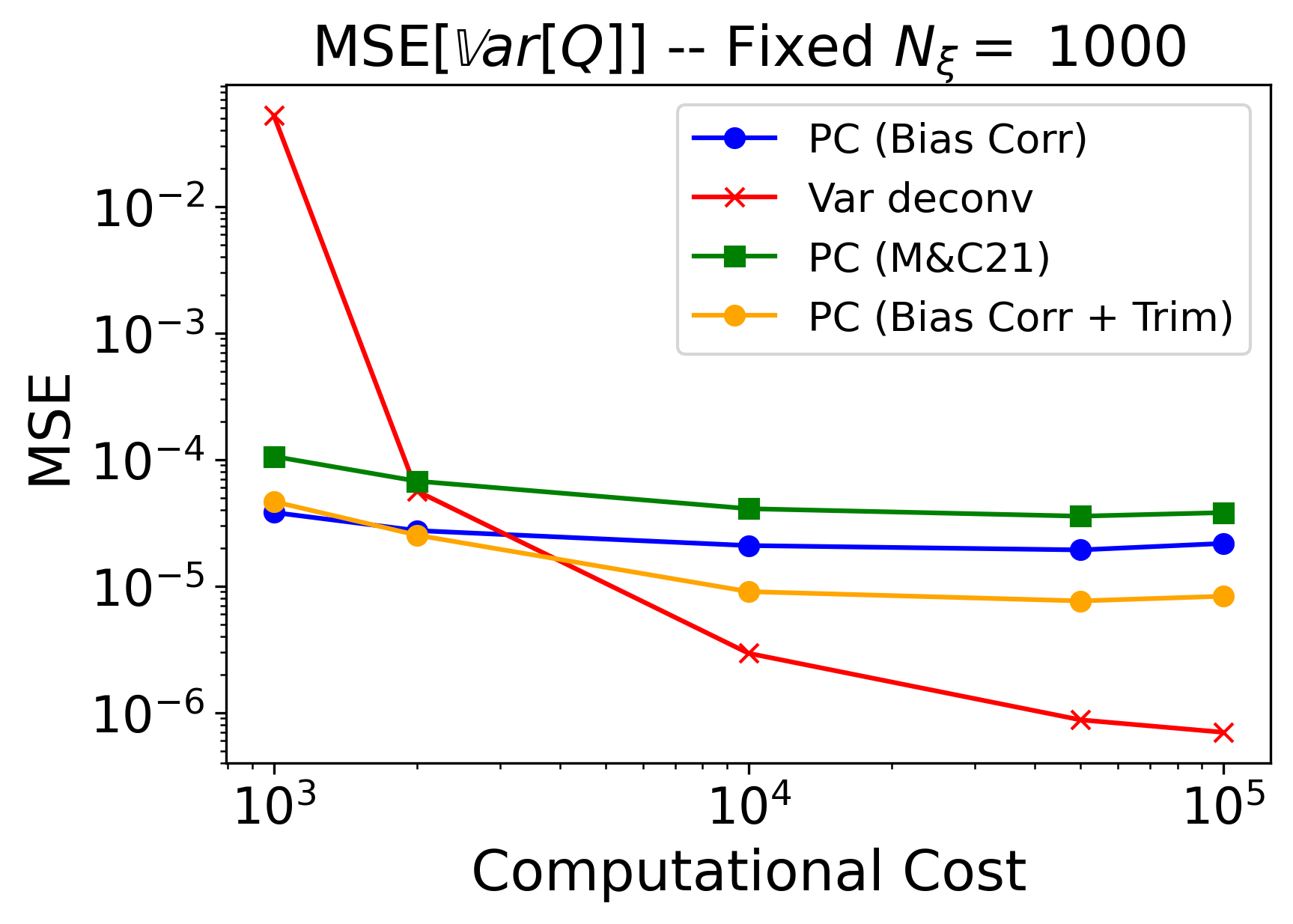}
         \caption{$\Nxi=1000$}
         \label{fig:MSE_var_Nxi1000}
     \end{subfigure}
     \hfill
     \begin{subfigure}[b]{0.32\textwidth}
         \centering
         \includegraphics[width=0.99\textwidth]{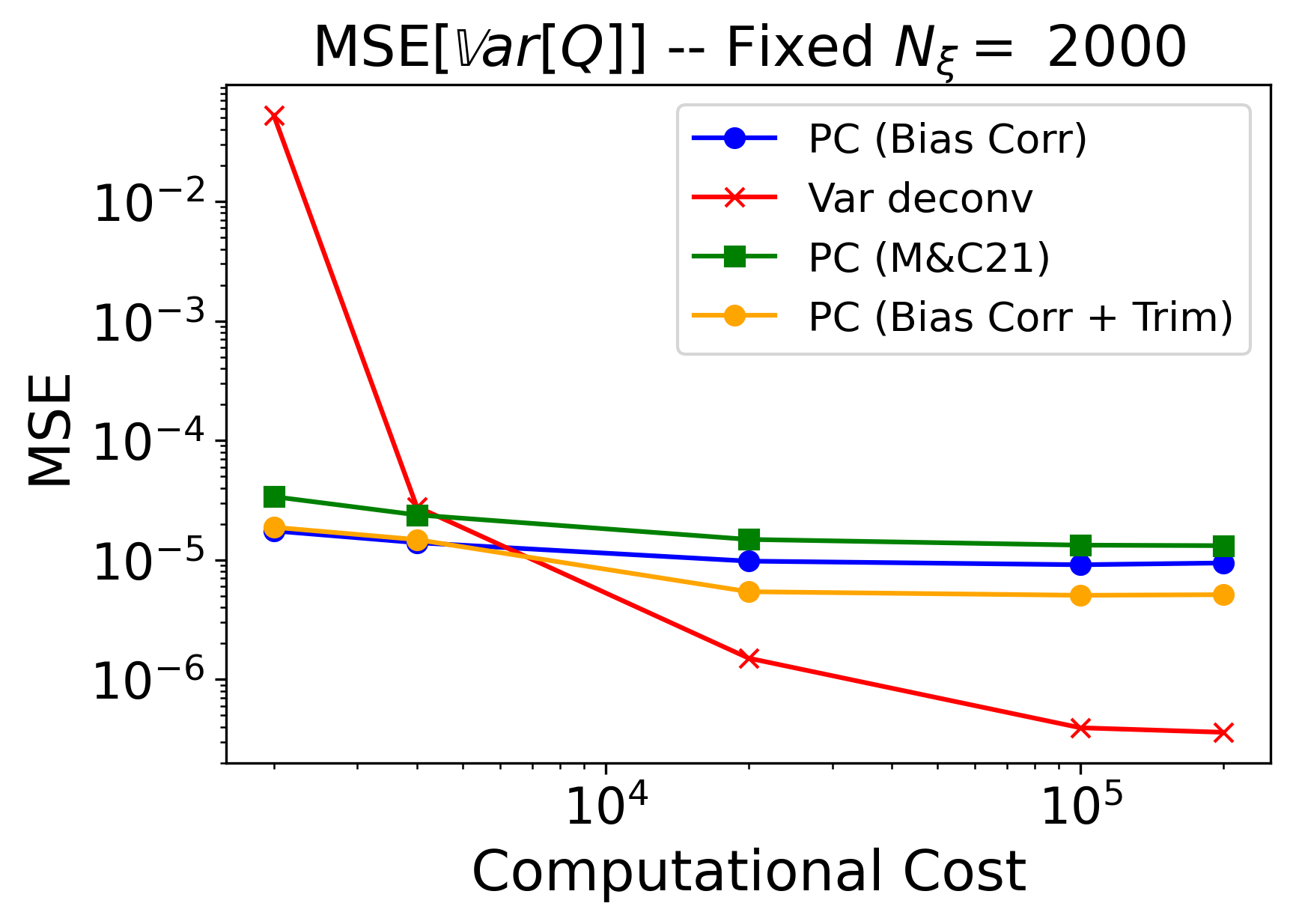}
         \caption{$\Nxi=2000$}
         \label{fig:MSE-var_nxi2000}
     \end{subfigure}
     \caption{MSE for the estimated variance obtained with $1\,500$ independent repetitions with an increasing number of UQ samples $\Nxi$ and $\Neta=\left[ 1, 2, 10, 50, 100\right]$.}
     \label{fig:MSE_var_Nxi}
\end{figure}

\begin{figure}[hbt!]
     \centering
     \begin{subfigure}[b]{0.32\textwidth}
         \centering
         \includegraphics[width=0.99\textwidth]{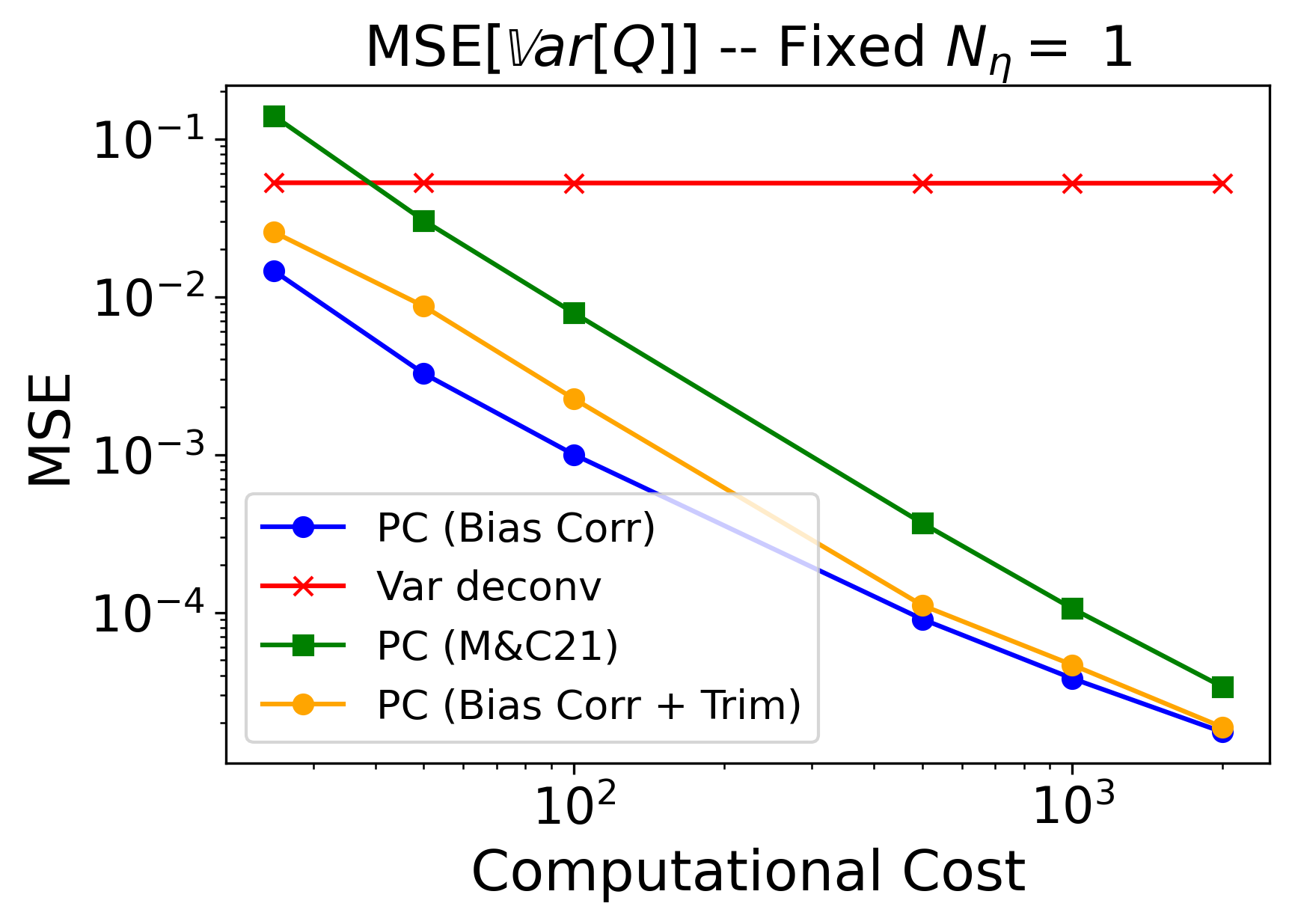}
         \caption{$\Neta=1$}
     \end{subfigure}
     \hfill
     \begin{subfigure}[b]{0.32\textwidth}
         \centering
         \includegraphics[width=0.99\textwidth]{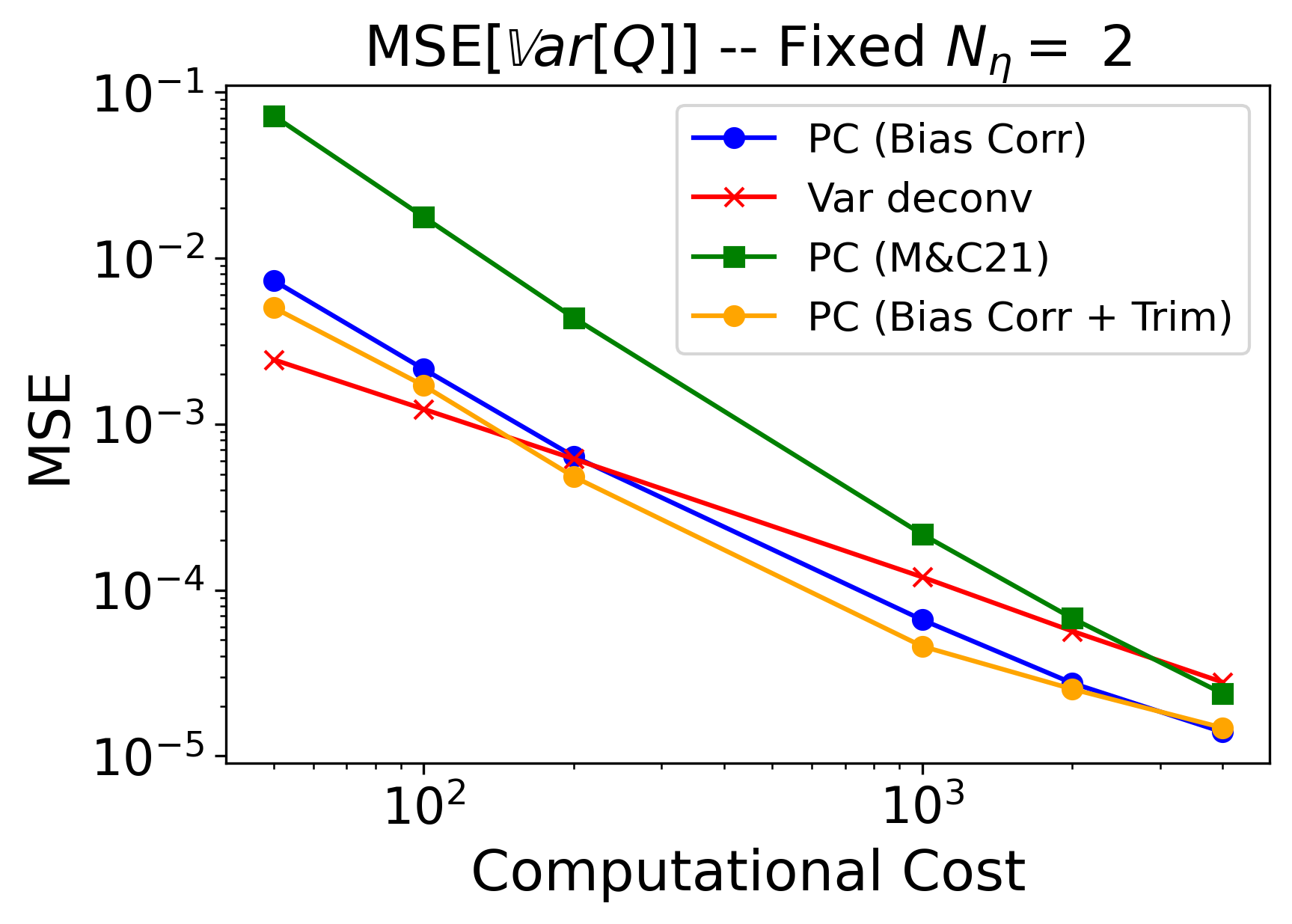}
         \caption{$\Neta=2$}
     \end{subfigure}
          \begin{subfigure}[b]{0.32\textwidth}
         \centering
         \includegraphics[width=0.99\textwidth]{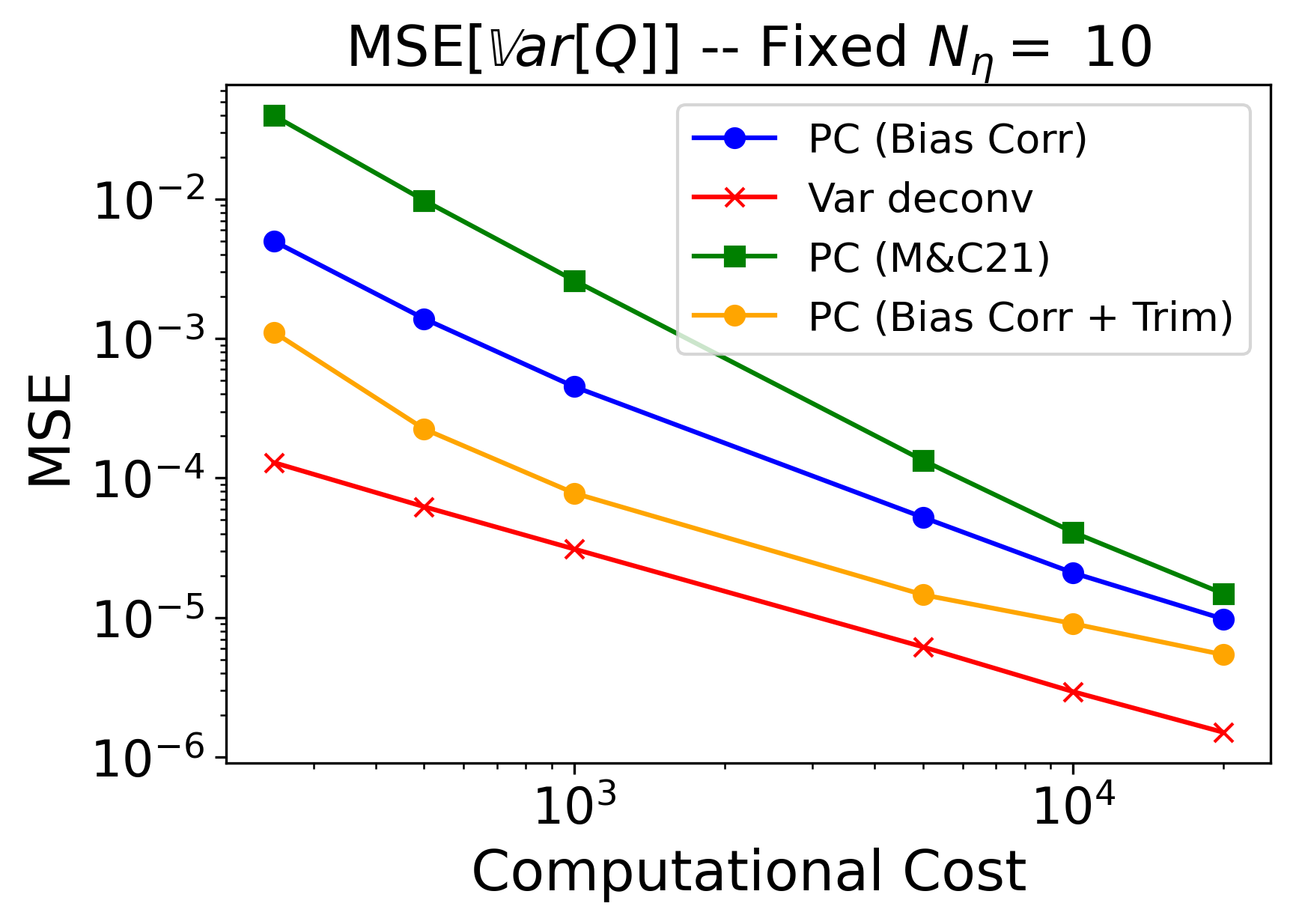}
         \caption{$\Neta=10$}
     \end{subfigure}
     \hfill
     \begin{subfigure}[b]{0.32\textwidth}
         \centering
         \includegraphics[width=0.99\textwidth]{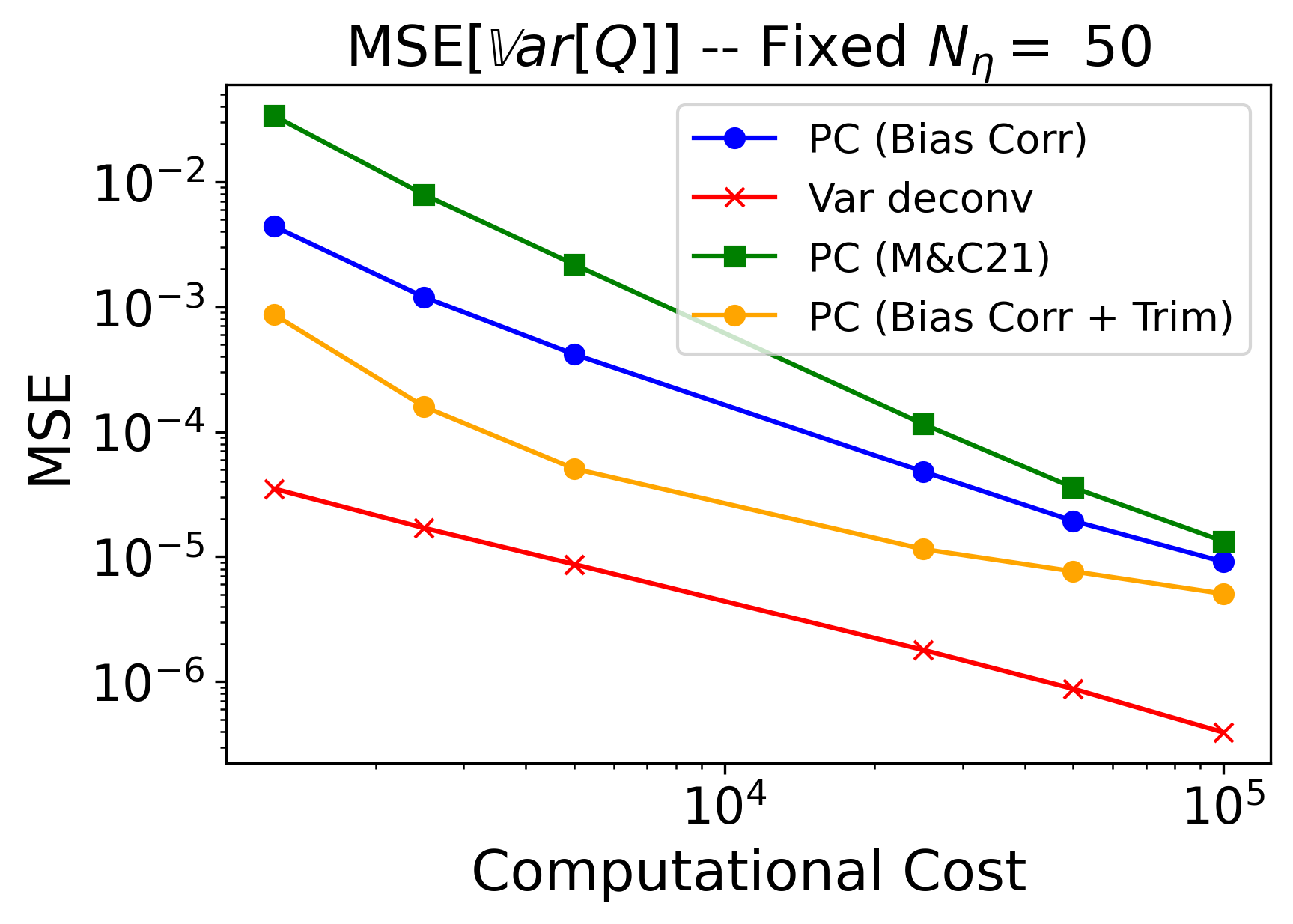}
         \caption{$\Neta=50$}
     \end{subfigure}
          \begin{subfigure}[b]{0.32\textwidth}
         \centering
         \includegraphics[width=0.99\textwidth]{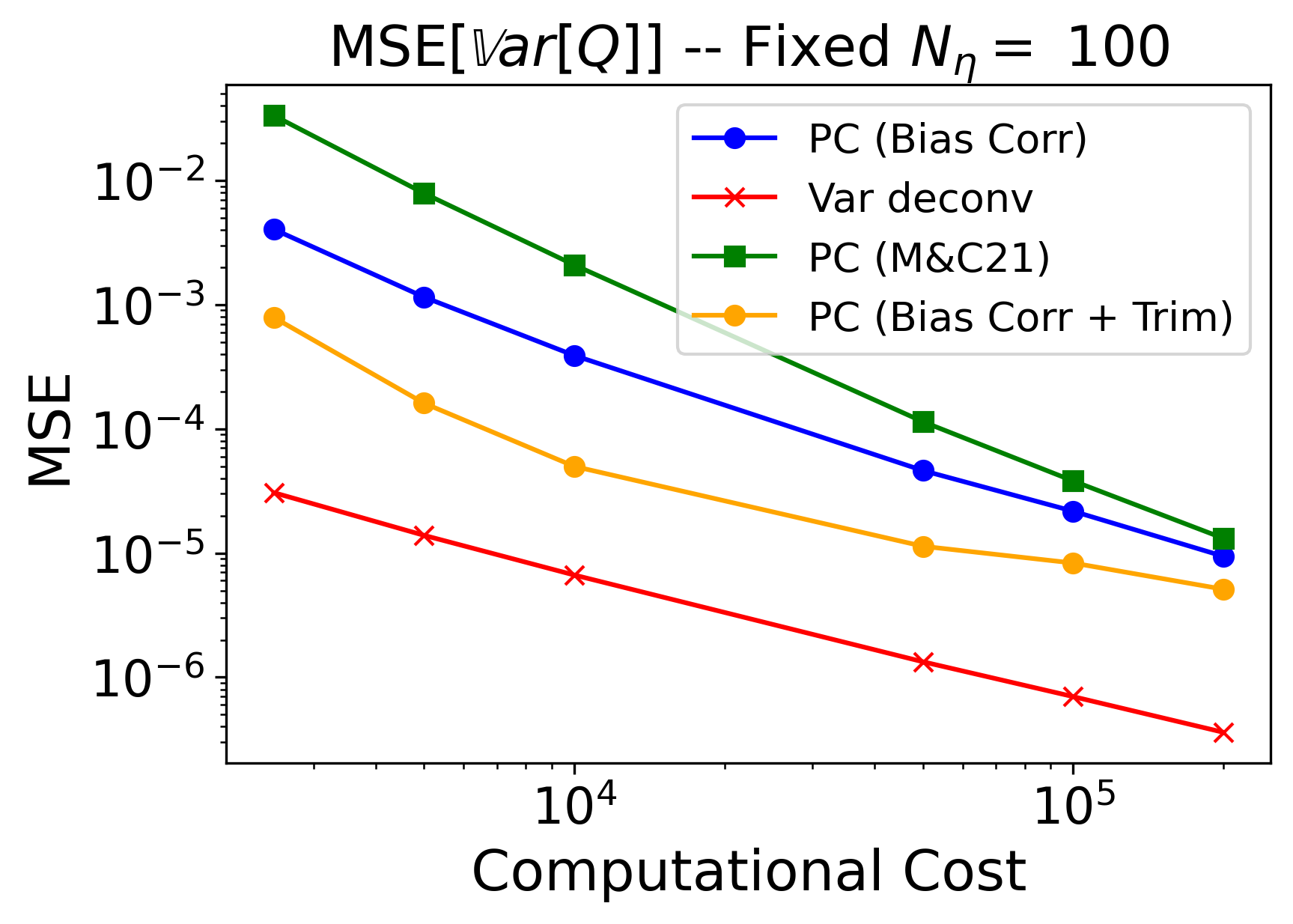}
         \caption{$\Neta=100$}
     \end{subfigure}
     \caption{MSE for the estimated variance obtained with $1\,500$ independent repetitions with an increasing number of particles $\Neta$ (per UQ sample) and $\Nxi=\left[ 25, 50, 100, 500, 1\,000, 2\,000\right]$.}
     \label{fig:MSE_var_Neta}
\end{figure}

\begin{figure}[hbt!]
     \centering
     \begin{subfigure}[b]{0.329\textwidth}
         \centering
         \includegraphics[width=0.99\textwidth]{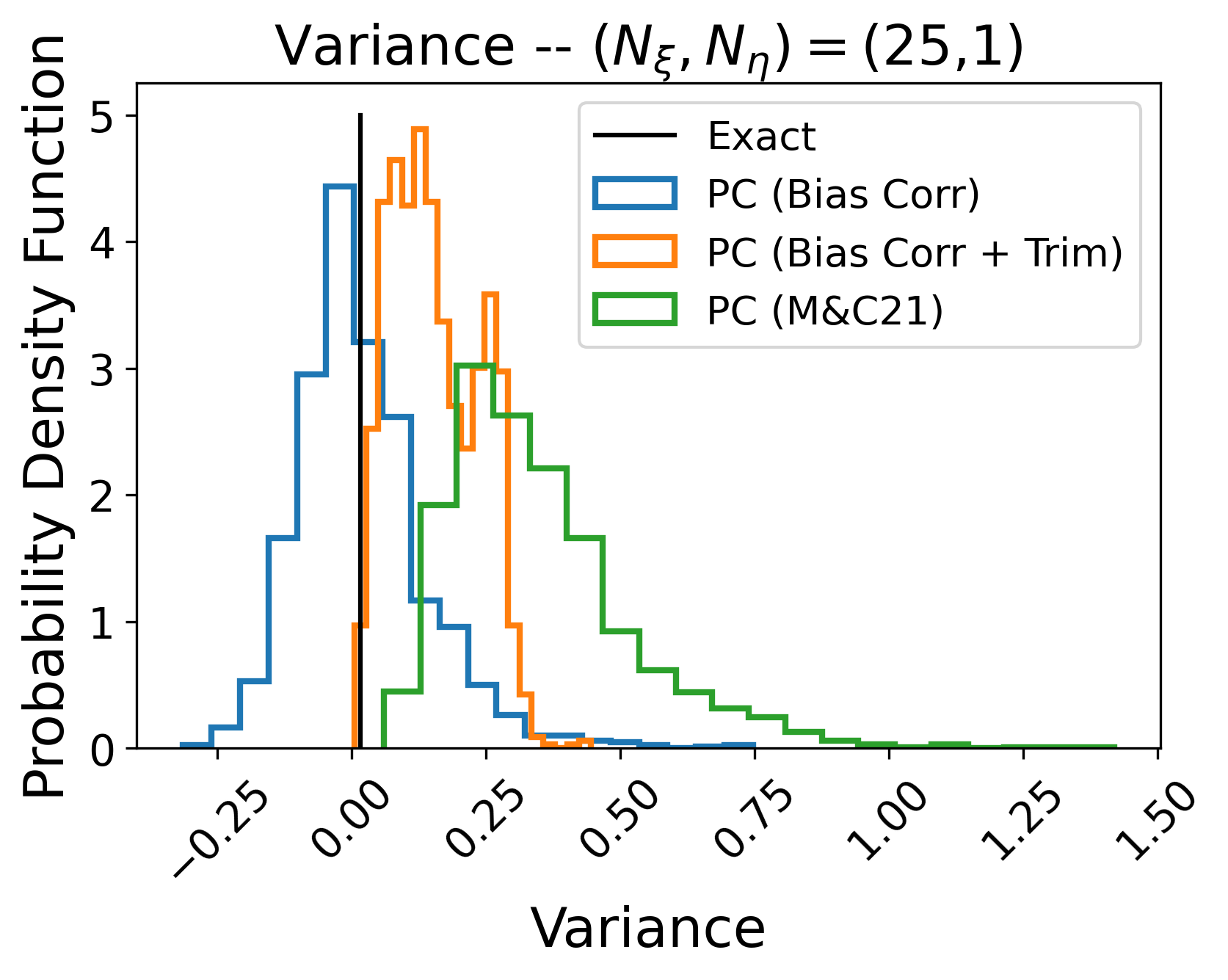}
         \vspace{-1mm}
         \caption{$(N_\xi,N_\eta)=(25,1)$}
         \label{fig:variance_all}
     \end{subfigure}
     \hfill
     \begin{subfigure}[b]{0.329\textwidth}
         \centering
         \includegraphics[width=0.99\textwidth]{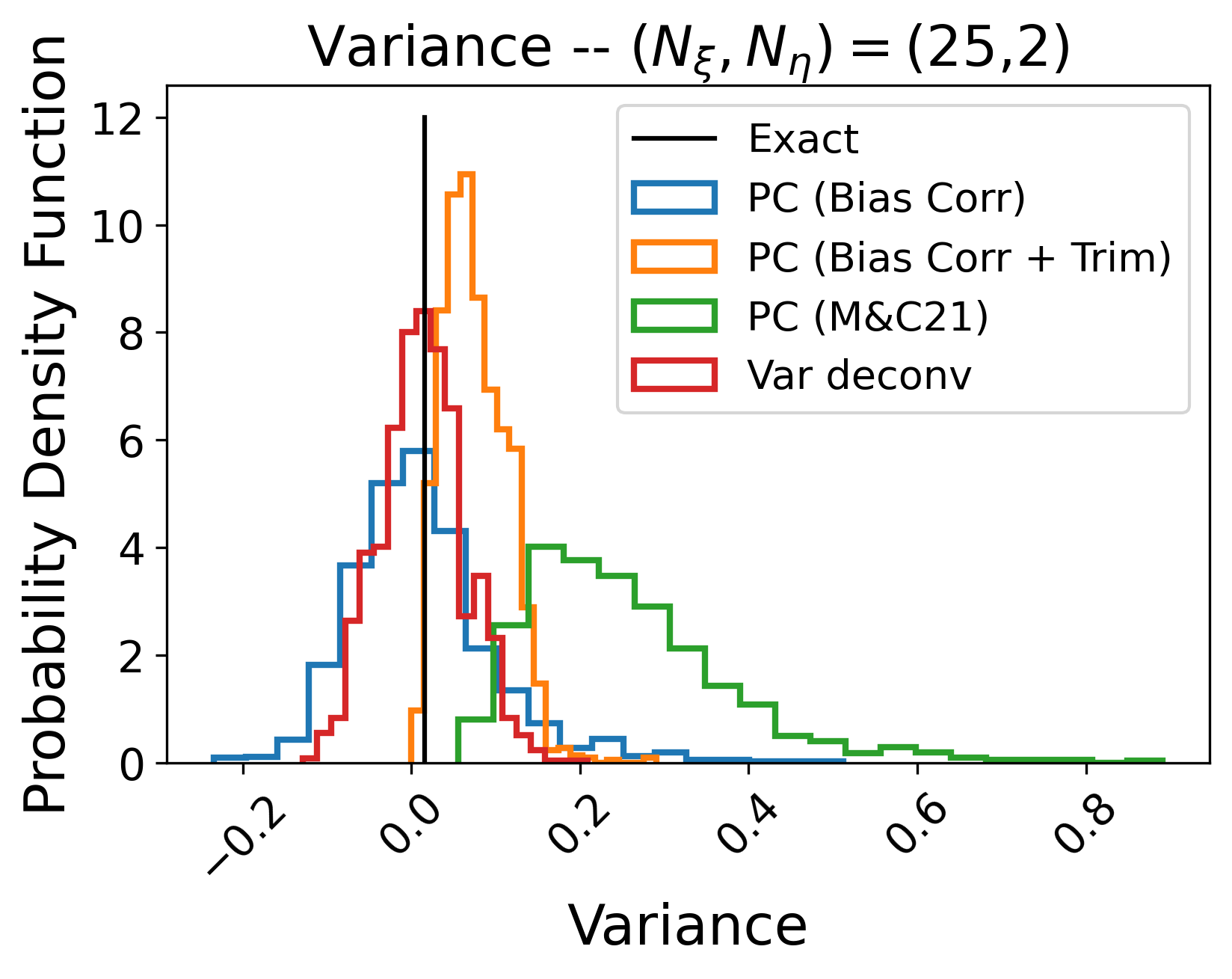}
         \vspace{-1mm}
         \caption{$(N_\xi,N_\eta)=(25,2)$}
         \label{fig:variance_all}
     \end{subfigure}
     \hfill
     \begin{subfigure}[b]{0.329\textwidth}
         \centering
         \includegraphics[width=0.99\textwidth]{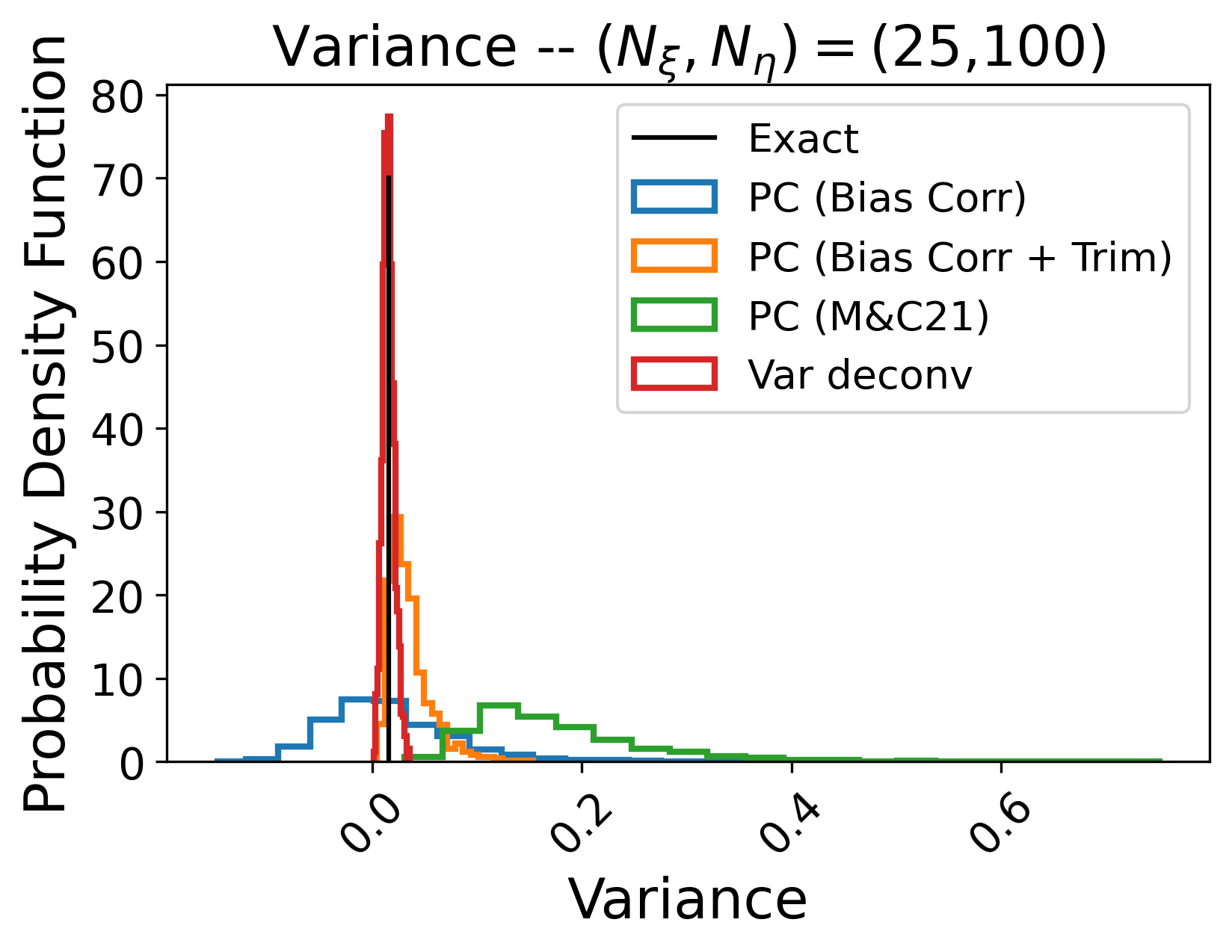}
         \vspace{-1mm}
         \caption{$(N_\xi,N_\eta)=(25,100)$}
         \label{fig:variance_all}
     \end{subfigure} 
     \begin{subfigure}[b]{0.329\textwidth}
     \vspace{1mm}
         \centering
         \includegraphics[width=0.99\textwidth]{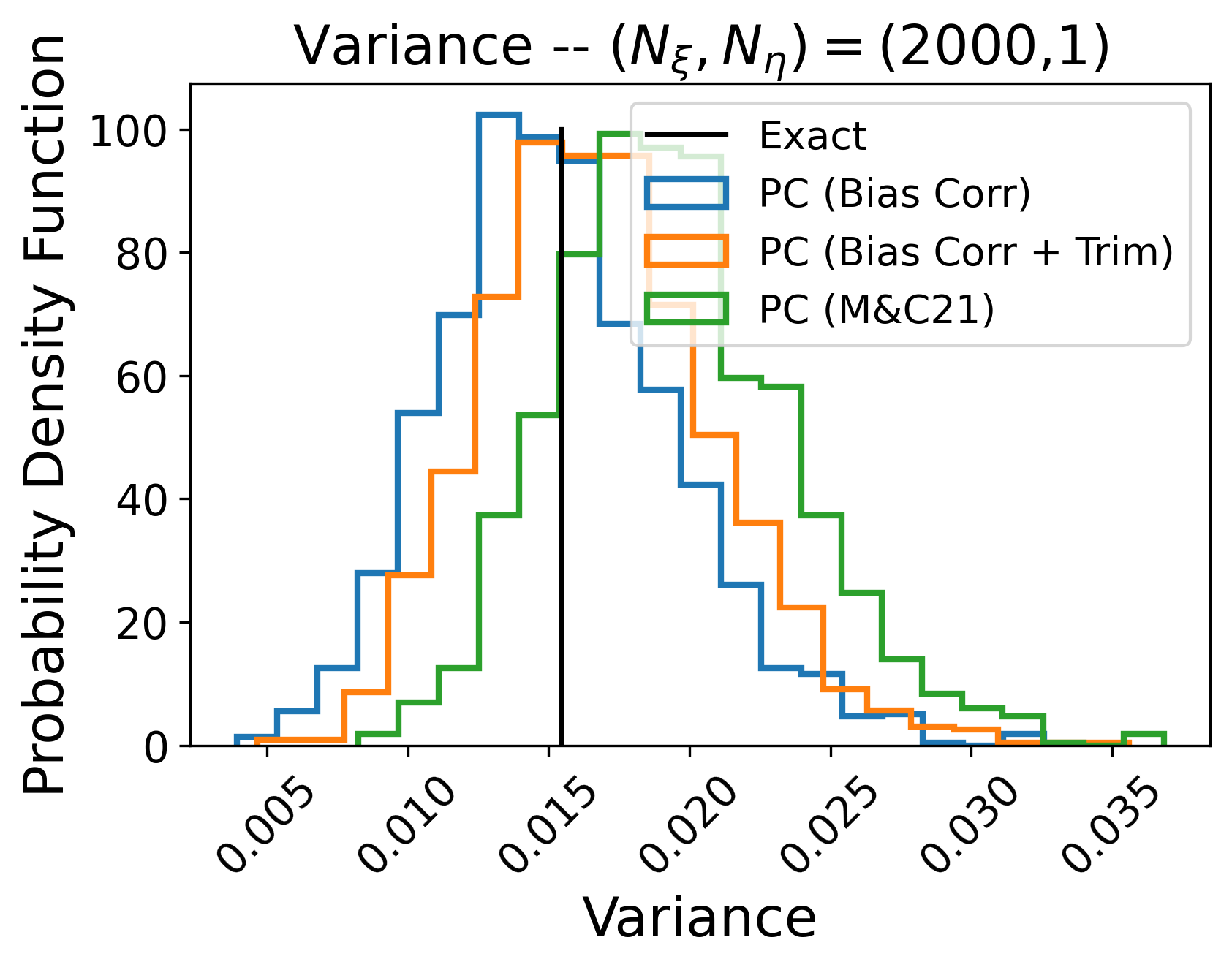}
         \vspace{-1mm}
         \caption{$(N_\xi,N_\eta)=(2000,1)$}
         \label{fig:variance_all}
     \end{subfigure}
     \hfill
     \begin{subfigure}[b]{0.329\textwidth}
     \vspace{1mm}
         \centering
         \includegraphics[width=0.99\textwidth]{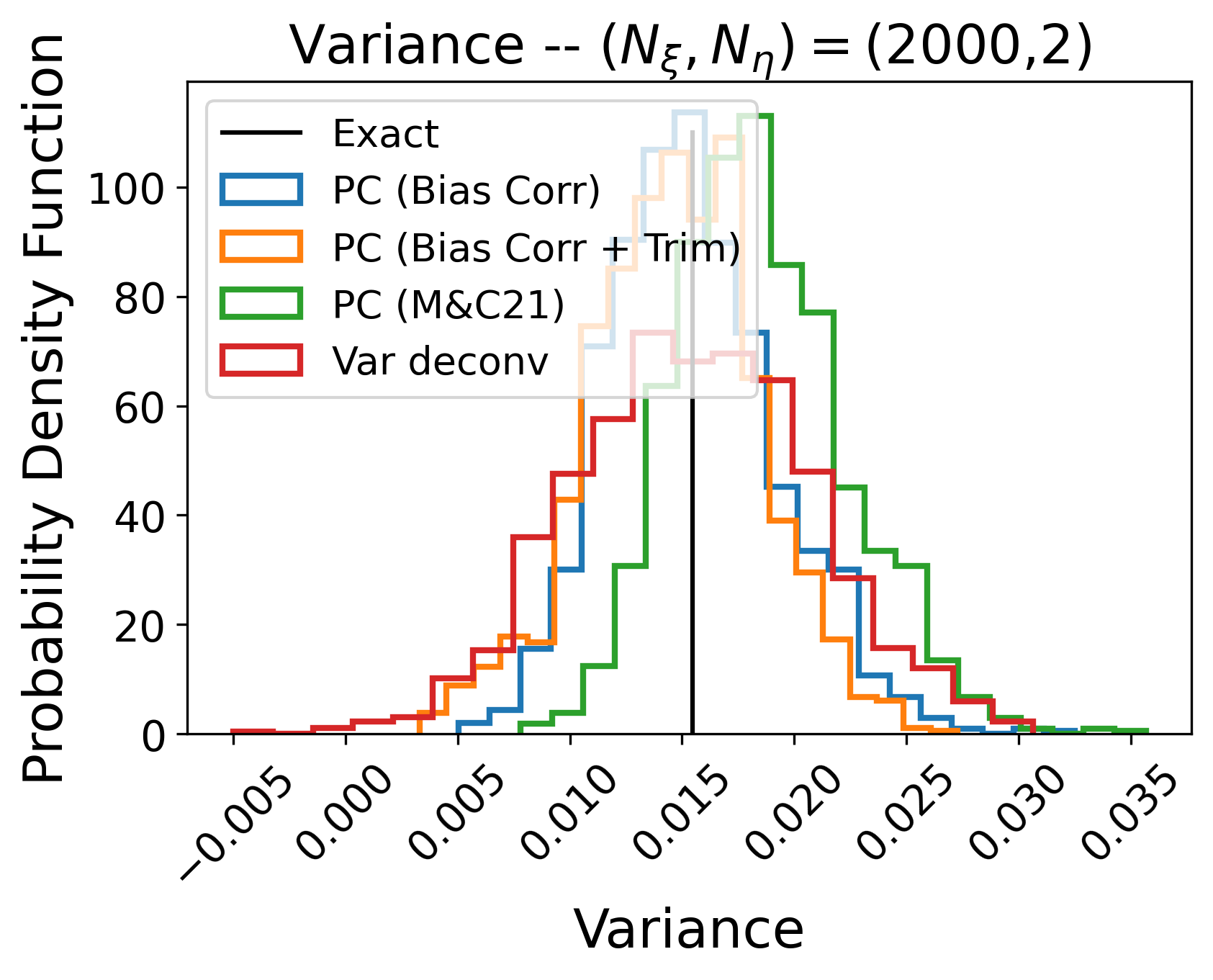}
         \vspace{-1mm}
         \caption{$(N_\xi,N_\eta)=(2000,2)$}
         \label{fig:variance_all}
     \end{subfigure}
     \hfill
     \begin{subfigure}[b]{0.329\textwidth}
     \vspace{1mm}
         \centering
         \includegraphics[width=0.99\textwidth]{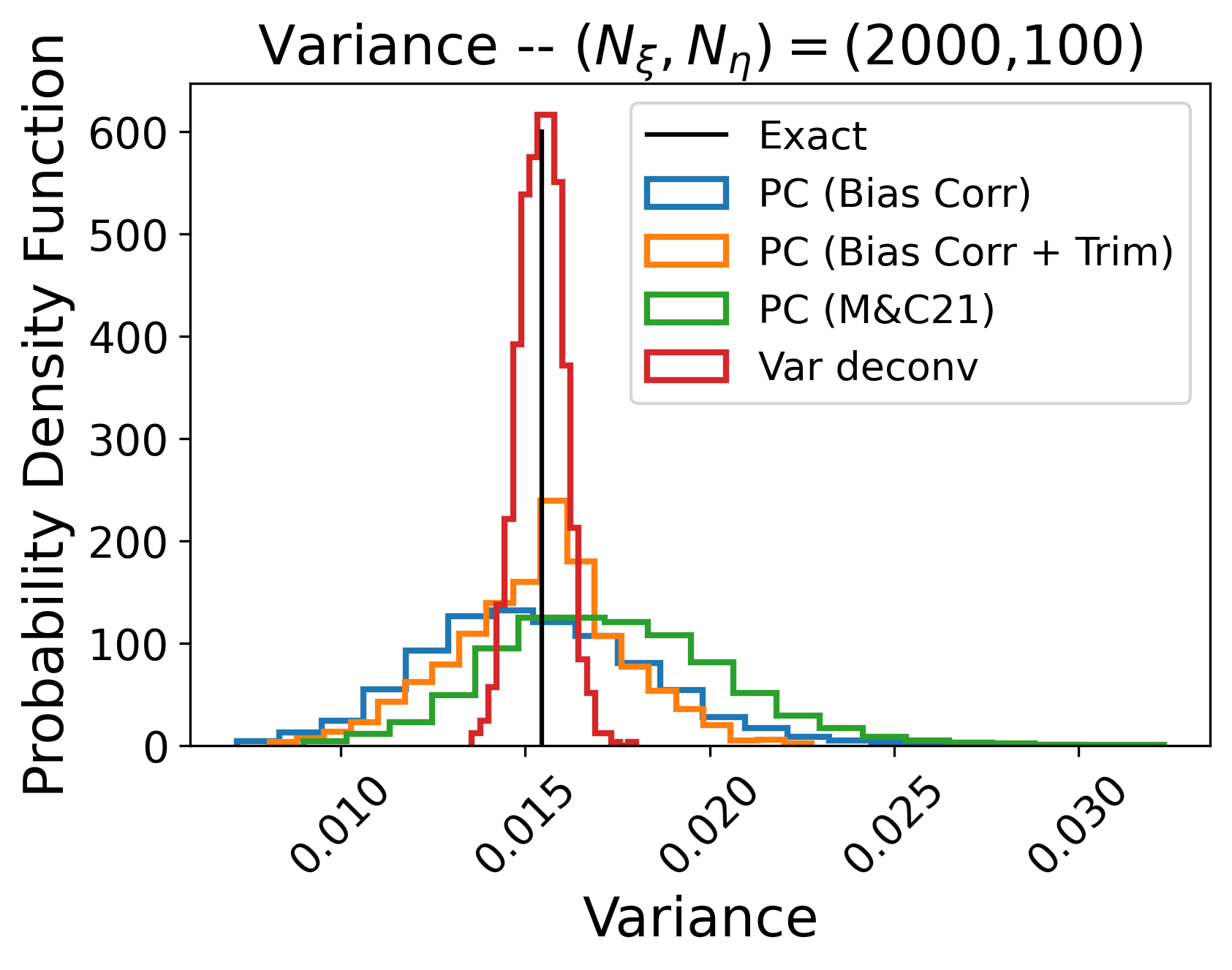}
         \vspace{-1mm}
         \caption{$(N_\xi,N_\eta)=(2000,100)$}
         \label{fig:variance_all}
     \end{subfigure}
     \caption{Probability density functions for the estimated variance with PCE and variance deconvolution.}
     \label{fig:pdf_var}
\end{figure}

We now move to the quantification of statistics, for the case with $d=3$, namely the variance of the transmittance for the problem with three uncertainties. We consider all the combinations of $(\Nxi,\Neta)$ given by the following set of numbers of UQ samples and particles, namely $\Nxi=\left[25, 50, 100, 500, 1\,000, 2\,000 \right]$ and $\Neta = \left[ 1, 2, 10, 50, 100 \right]$. We consider the PC method as presented in~\cite{GeraciMC2021} (\emph{PC M\&C21}), the PC method with the addition of the bias correction introduced in Sec.~\ref{SEC:bias_corr} (\emph{PC (Bias Corr)}) and with the addition of bias and expansion trim (\emph{PC (Bias Corr + Trim)}). Moreover, to provide a comparison with a sampling scheme, we also report the value computed using variance deconvolution method (\emph{Var deconv}) presented in~\cite{ClementsANS2022}. In Figs.~\ref{fig:MSE_var_Nxi} and~\ref{fig:MSE_var_Neta} we report all the results with fixed $\Nxi$ or $\Neta$. For all these results, it is possible to note that the algorithmic improvements introduced in this contribution significantly improve on the results of the original approach introduced in~\cite{GeraciMC2021}. Moreover, it is important to note how the PC is able to provide significantly better results than the deconvolution methods whenever the MC RT runs are severely under-resolved, i.e., with $\Neta=1$ or $2$. Otherwise, the sampling approach is more effective. Notably, the objective of this paper is to demonstrate an improvement with respect to~\cite{GeraciMC2021}, while a more rigorous comparison with different methods is left for a future paper. Indeed, the comparison among approaches would require a large array of test cases with different features since distinct methods could be best suited in specific circumstances, e.g., as a function of the dimension and regularity of the problem.
To further highlight the difference among methods, we report the probability density functions obtained for the $1\,500$ repetitions for the lowest and greater number of UQ samples and for the lowest and largest number of particles, plus $\Neta=2$, which is the minimum number of particles that allows for the use of variance deconvolution. The results are reported in Fig.~\ref{fig:pdf_var}. Two observations are in order after observing Fig.~\ref{fig:pdf_var}. First, all the methods exhibit distributions with negligible bias for the larger number of UQ samples, i.e., $\Nxi=2000$. Second, whenever the number of UQ samples is coarse, i.e., $\Nxi=25$, we observe a bias in the PC results of~\cite{GeraciMC2021}, whereas the introduction of unbiased estimation of the squared coefficients eliminates the bias, as predicted by the above theory. A biased estimator is obtained by using the expansion trim because not all the $1\,500$ independent realizations contain the exact same set of coefficients, due to the statistical variability. However, the introduction of the expansion trim ensures that the variance is always positive; this latter feature is not guaranteed in the bias-corrected approach because non-biased estimators of small quantities close to zero cannot be made always positive.  

\subsection{Global Sensitivity Analysis}
\begin{figure}
     \centering
     \begin{subfigure}[b]{0.49\textwidth}
         \centering
         \includegraphics[width=0.79\textwidth]{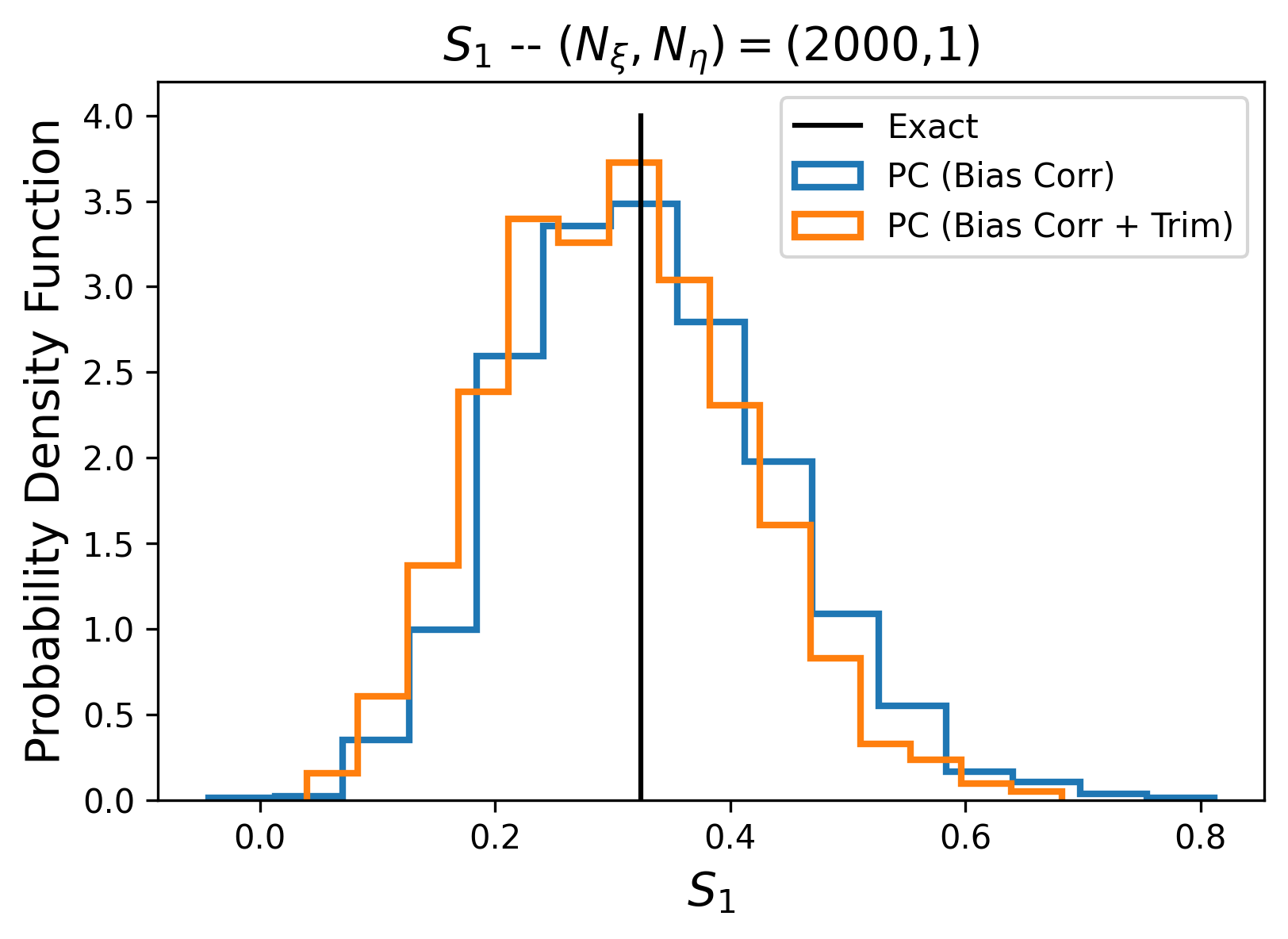}
         \vspace{-3mm}
         \caption{Sensitivity index $S_1$}
         \label{fig:S1}
     \end{subfigure}
     \hfill
     \begin{subfigure}[b]{0.49\textwidth}
         \centering
         \includegraphics[width=0.79\textwidth]{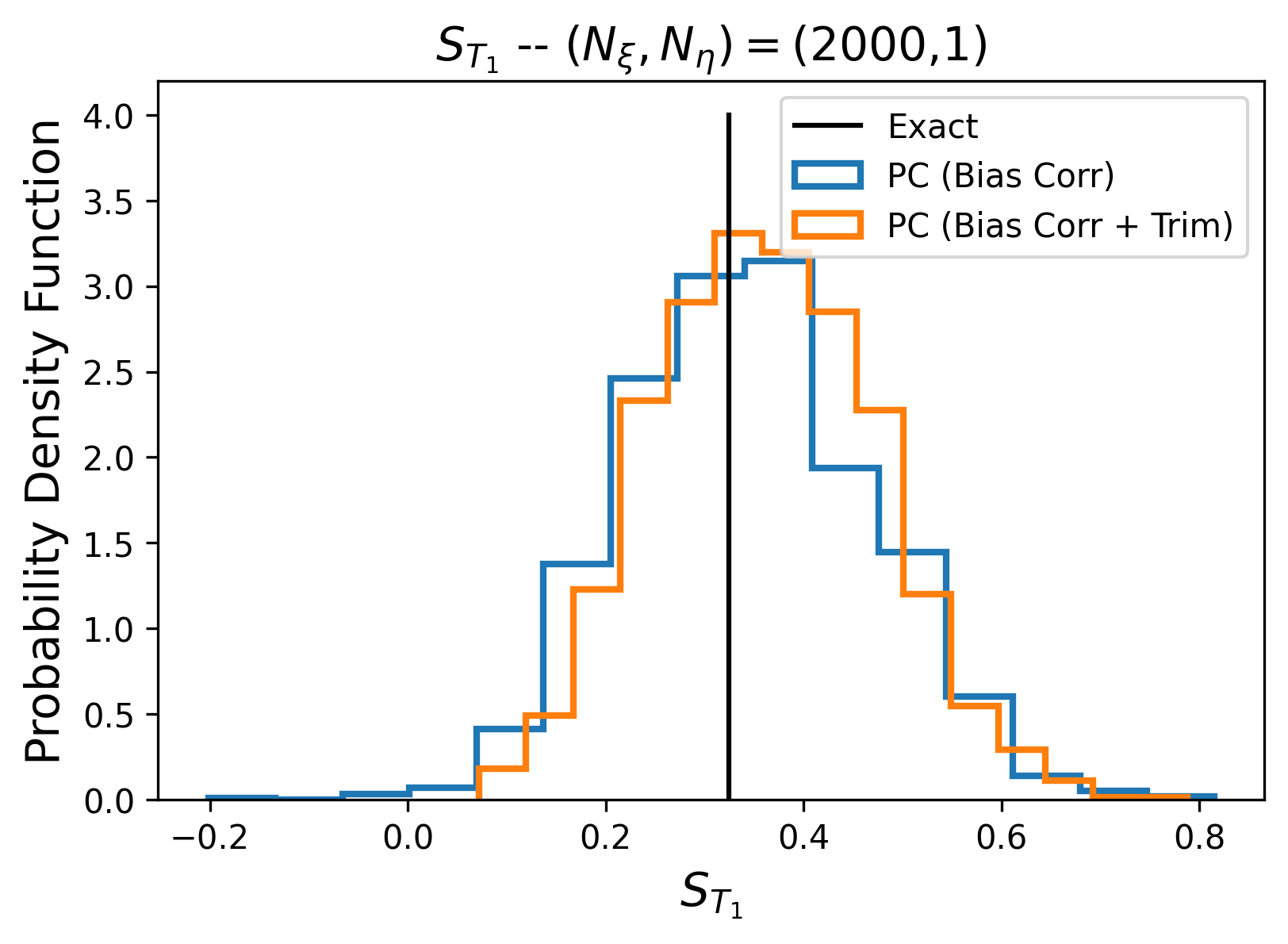}
         \vspace{-3mm}
         \caption{Total sensitivity index $S_{T_1}$}
         \label{fig:TS1}
     \end{subfigure}     
     \caption{Sensitivity index $S_1$~(\subref{fig:S1}) and total sensitivity index $S_{T_1}$ (\subref{fig:TS1}) obtained with the PC with bias correction and bias correction and expansion trim.}
     \label{fig:PCE_GSA}
\end{figure}

For completeness, we report here the GSA results obtained for the sensitivity indices associated to the first uncertain variable, i.e., the total cross section of the first slab section. The results are reported in Fig.~\ref{fig:PCE_GSA} and correspond to $1\,500$ independent realizations with $\Nxi=2\,000$ and $\Neta=1$, for the algorithm with bias correction and with and without the expansion trim.  

\section{CONCLUSIONS} 
\label{SEC:conclusions}
In this contribution, we illustrated several improvements with respect to the PC algorithm introduced in~\cite{GeraciMC2021}. We have discussed the efficiency of the PC coefficients' estimators in the presence of re-sampling cost, and we provided a bias correction expression for the expansion's terms in the variance. Moreover, we demonstrated how the PC variance can be used to estimate the PC variability due to the UQ random sampling, while correcting for the finite number of particles $\Neta$. Finally, we reported GSA results to illustrate a possible direction for future investments. We note here that this manuscript is part of a larger research activity in which several methods are being developed and compared to serve under different scenarios. For instance, a sampling method based on variance deconvolution for GSA is provided in the companion paper~\cite{ClementsMC2023}. A large array of test cases will be required to gauge the trade-off among different approaches, including some not explored here like quadrature- and regression-based approaches; this task is left for future studies. Finally, we note that the formulation presented here is amenable for the multilevel extension presented in~\cite{Merritt2020}, which aims at accelerating the PC construction whenever limited evaluations of a system are available, but a larger set of samples can be readily, or more efficiently, obtained from coarser low-fidelity models.

\section*{ACKNOWLEDGMENTS}
This article has been authored by an employee of National Technology \& Engineering Solutions of Sandia, LLC under Contract No. DE-NA0003525 with the U.S. Department of Energy (DOE). The employee owns all right, title and interest in and to the article and is solely responsible for its contents. The United States Government retains and the publisher, by accepting the article for publication, acknowledges that the United States Government retains a non-exclusive, paid-up, irrevocable, world-wide license to publish or reproduce the published form of this article or allow others to do so, for United States Government purposes. The DOE will provide public access to these results of federally sponsored research in accordance with the DOE Public Access Plan \texttt{https://www.energy.gov/downloads/doe-public-access-plan}.

\newif\ifusebibtex
\usebibtextrue

\ifusebibtex
\setlength{\baselineskip}{12pt}
\bibliographystyle{mc2023}
\bibliography{mc2023}
\else
\setlength{\baselineskip}{12pt}

\fi


\end{document}